\documentclass[a4paper,11pt]{article}
\usepackage[utf8]{inputenc}
\usepackage{a4wide}
\usepackage{latexsym,amsfonts,amsmath,amssymb,mathrsfs,url,amsthm}
\usepackage{mathtools}
\usepackage{color,graphicx}
\usepackage{comment}
\usepackage[subrefformat=parens]{subcaption}
\usepackage[toc]{appendix}
\providecommand{\keywords}[1]
{
  \noindent \small	
  \textbf{Keywords:} #1
}
\providecommand{\amscode}[1]
{
  \noindent \small	
  \textbf{AMS subject classifications:} #1
}

\newtheorem{theorem}{Theorem}
\newtheorem{lemma}{Lemma}
\newtheorem{example}{Example}

\newtheorem{remark}{Remark}

\newtheorem{corollary}{Corollary}

\newcommand{\EE}{\mathbb{E}}
\newcommand{\NN}{\mathbb{N}}
\newcommand{\RR}{\mathbb{R}}
\newcommand{\VV}{\mathbb{V}}
\newcommand{\bszero}{\boldsymbol{0}}
\newcommand{\bsone}{\boldsymbol{1}}
\newcommand{\Ncal}{\mathcal{N}}
\newcommand{\Xcal}{\mathcal{X}}

\allowdisplaybreaks

\usepackage{xcolor}
\newcommand{\tg}[1]{{\color{black}{#1}}} 
\newcommand{\tgg}[1]{{\color{black}{#1}}}

\title{Constructing unbiased gradient estimators with finite variance for conditional stochastic optimization\thanks{The work of the first author was supported by JSPS KAKENHI Grant Number 20K03744.}}
\author{Takashi Goda\thanks{School of Engineering, University of Tokyo, 7-3-1 Hongo, Bunkyo-ku, Tokyo 113-8656, Japan ({\tt goda@frcer.t.u-tokyo.ac.jp}; {\tt wataru.kitade@gmail.com})}, Wataru Kitade\footnotemark[2]}
\date{\today}

\begin{document}

\maketitle

\begin{abstract}
    We study stochastic gradient descent for solving conditional stochastic optimization problems, in which an objective to be minimized is given by a parametric nested expectation with an outer expectation taken with respect to one random variable and an inner conditional expectation with respect to the other random variable. The gradient of such a parametric nested expectation is again expressed as a nested expectation, which makes it hard for the standard nested Monte Carlo estimator to be unbiased. In this paper, we show under some conditions that a multilevel Monte Carlo gradient estimator is unbiased and has finite variance and finite expected computational cost, so that the standard theory from stochastic optimization for a parametric (non-nested) expectation directly applies. We also discuss a special case for which yet another unbiased gradient estimator with finite variance and cost can be constructed.
\end{abstract}

\amscode{65C05, 65K10, 90C15}

\section{Introduction}\label{sec:intro}
In this paper we consider optimization problems of the form
\begin{align}\label{eq:CSO} 
    \min_{x\in \Xcal}F(x)\quad \text{with}\quad F(x)=\EE_{\xi}\left[ f_{\xi}\left( \EE_{\eta\mid \xi}\left[ g_{\eta}\left(x, \xi\right)\right]\right)\right],
\end{align}
where $\Xcal\subseteq \RR^d$ denotes a feasible solution region, $\xi$ and $\eta$ are (not necessarily independent) random variables, $g_{\eta}(\cdot,\xi): \RR^d\to \RR^k$ is a vector-valued differentiable function dependent both on $\xi$ and $\eta$, and $f_{\xi}(\cdot): \RR^k\to \RR$ is a real-valued differentiable function dependent on $\xi$ but not on $\eta$. This type of problems is called \emph{conditional stochastic optimization (CSO)} \cite{HCH20,HZCH20,HCH21} and appear in a variety of different fields of application, such as instrumental variable (IV) regression \cite{HLLT17,SSG19,MMLR20}, model-agnostic meta-learning \cite{FAL17,HZCH20}, and Bayesian experimental design \cite{CV95,BDET20,CDRLT20,GHI20,GHKF22}. \tg{As one of the motivating examples, let us describe a non-linear IV regression problem \cite{BKS19,MMLR20}.
\begin{example}\label{exm:iv}
Consider the following data generating process
\[ Y=f(X)+\varepsilon,\]
with an unknown, non-linear continuous function $f$ and the additive noise $\varepsilon$ satisfying $\EE[\varepsilon]=0$ and $\EE[\varepsilon^2]<\infty$. The noise $\varepsilon$ is correlated with the input variable $X$ so that $\EE[\varepsilon\mid X]\neq 0$ in general. Let us assume that there exists an accessible instrumental variable $Z$ which satisfies the conditions that (1) $Z$ has a causal influence on $X$, that (2) $Z$ affects $Y$ only through $X$ and that (3) $Z$ is independent of the noise $\varepsilon$, which implies $\EE[\varepsilon\mid Z]= 0$. The problem of our interest is then to estimate $f$, for instance, by using a dataset $\{(X_i,Y_i,Z_i)\}_{i=1,\ldots,n}$. In \cite{MMLR20}, this IV problem is formulated as
\begin{align}\label{eq:CSO_IV}
    \min_{f}R(f)\quad \text{with}\quad R(f)=\EE_{Y,Z}\left[ L\left(Y,\EE_{X\mid Z}\left[ f(X)\right]\right)\right],
\end{align} 
where $L(x,y)$ measures a difference between $x$ and $y$, e.g., $L(x,y)=\|x-y\|_2^2$. If we model a function $f$ by a set of parameters $\theta$, this optimization problem can be regarded as a concrete example of \eqref{eq:CSO}.
\end{example}}

The essential difficulty of solving CSO problems comes from the fact that the objective function $F$ is defined as a nested expectation with an outer expectation with respect to $\xi$ and an inner conditional expectation with respect to $\eta$. In fact, the chain rule applies under mild conditions and we see that the gradient of $F$ with respect to $x$ is again given by a nested expectation
\[ \nabla F(x) =  \EE_{\xi}\left[ \left( \EE_{\eta\mid \xi}\left[ \nabla g_{\eta}\left(x, \xi\right)\right]\right)^{\top}\nabla f_{\xi}\left( \EE_{\eta\mid \xi}\left[ g_{\eta}\left(x, \xi\right)\right]\right)\right], \]
so that it is not clear how to construct an efficient gradient estimator for stochastic gradient descent, which is not only unbiased but also has finite variance under finite computational cost. Constructing biased gradient estimators are relatively straightforward by ``nesting" the naive Monte Carlo estimators for the inner and outer expectations \cite{RCYWW18}:
\begin{align}\label{eq:NMC}
    \nabla F(x) \approx \widehat{\nabla F}_{M,N}(x) = \frac{1}{N}\sum_{n=1}^{N} \left( \frac{1}{M}\sum_{m=1}^{M} \nabla g_{\eta^{(m,n)}}(x, \xi^{(n)})\right)^{\top} \nabla f_{\xi}\left( \frac{1}{M}\sum_{m=1}^{M} g_{\eta^{(m,n)}}(x, \xi^{(n)})\right),
\end{align}
with $\xi^{(1)},\ldots,\xi^{(N)}$ being i.i.d.\ copies of $\xi$, and for each $n=1,\ldots,N$, $\eta^{(1,n)},\ldots,\xi^{(M,n)}$ being i.i.d.\ copies of $\eta$ conditioned on $\xi=\xi^{(n)}$. Using such a biased gradient estimator in stochastic gradient descent leads to a bias in the solution and we need to be quite careful about the bias-variance trade-off in gradient estimation to establish a theoretical guarantee of convergence to a globally optimal solution, see \cite{HZCH20}. The simple test case given in \cite[Section~4.1]{GHKF22} clearly depicts a failure of the nested Monte Carlo gradient estimator with fixed small $M$ both in theory and computation.

To tackle this issue of CSO problems, there have been several attempts to apply (suitably randomized) multilevel Monte Carlo (MLMC) methods to construct de-biased gradient estimators, such as \cite{BGILZ17,HCH21} in a general CSO setting as well as \cite{GHKF22,IG21,SC21} in more specific contexts, respectively. MLMC methods were originally introduced by Heinrich in \cite{H98} for parametric integration and then by Giles in \cite{G08} for estimating the expectations arising from stochastic differential equations. A randomization technique to make MLMC estimators unbiased has been proposed \tg{independently in \cite{M11} and \cite{RG15}.} Besides CSO problems, MLMC methods have been applied \tg{in the contexts of stochastic approximation in \cite{F16,BGP19,DG19} and filtering problems in \cite{JKLZ17,JKOZ18,JLY22}, respectively,} and also used to estimate nested expectations efficiently in \cite{BHR15,G15,GG19,GH19,GHI20,HGGT20} as well as some others. Moreover, constructing unbiased MLMC gradient estimators has been studied quite recently in \cite{JLL21,BJ22} for a class of problems different from what we consider in this paper.

The aim of this paper is to push forward the theoretical analysis made in \cite{HCH21} further and to show under some conditions that the gradient $\nabla F(x)$ can be estimated without any bias and with finite variance and finite expected computational cost. Our result enables the standard theory from stochastic optimization to apply directly without any need to be cautious about a delicate bias-variance trade-off in computation. We also discuss a special case of CSO problems for which yet another unbiased gradient estimator can be constructed.
\section{MLMC gradient estimation}\label{sec:mlmc}
\subsection{Infinite series representation}
Let us consider the nested Monte Carlo gradient estimator with $N=1$ and $M=2^{\ell}$ for some $\ell\in \NN\cup\{0\}$ and denote it by $\psi_\ell$, that is,
\[ \psi_{\ell}(x):=\widehat{\nabla F}_{2^{\ell},1}(x)=\left( \frac{1}{2^{\ell}}\sum_{m=1}^{2^{\ell}}\nabla g_{\eta^{(m)}}(x,\xi)\right)^{\top}\nabla f_{\xi}\left( \frac{1}{2^{\ell}}\sum_{m=1}^{2^{\ell}}g_{\eta^{(m)}}(x,\xi)\right), \]
with $\eta^{(1)},\ldots,\eta^{(2^{\ell})}$ being i.i.d.\ copies of $\eta$ conditioned on $\xi$. For the sake of notational simplicity, we write 
\[ \overline{g_{\bullet}(x,\xi)}^{\ell}=\frac{1}{2^{\ell}}\sum_{m=1}^{2^{\ell}}g_{\eta^{(m)}}(x,\xi), \quad \overline{\nabla g_{\bullet}(x,\xi)}^{\ell}=\frac{1}{2^{\ell}}\sum_{m=1}^{2^{\ell}}\nabla g_{\eta^{(m)}}(x,\xi), \]
to represent the respective inner Monte Carlo averages based on $2^{\ell}$ i.i.d.\ copies of $\eta$. Then we have
\[ \psi_{\ell}(x)=\left(\overline{\nabla g_{\bullet}(x,\xi)}^{\ell}\right)^{\top}\nabla f_{\xi}\left(\overline{g_{\bullet}(x,\xi)}^{\ell}\right).\]
Under some mild measurability conditions, \tg{Lebesgue’s dominated convergence theorem} ensures that
\[ \nabla F(x) = \lim_{\ell\to \infty}\EE\left[ \psi_{\ell}(x)\right]  = \lim_{L\to \infty}\left(\EE\left[ \psi_{0}\right] + \sum_{\ell=1}^{L}\EE\left[ \psi_{\ell}(x)-\psi_{\ell-1}(x)\right]\right).\]

For the same variables $\eta^{(1)},\ldots,\eta^{(2^{\ell})}$ used in $\psi_\ell$, we write
\begin{align*}
& \overline{g_{\bullet}(x,\xi)}^{\ell-1,(a)} = \frac{1}{2^{\ell-1}}\sum_{m=1}^{2^{\ell-1}}g_{\eta^{(m)}}(x,\xi), \quad \overline{\nabla g_{\bullet}(x,\xi)}^{\ell-1,(a)} = \frac{1}{2^{\ell-1}}\sum_{m=1}^{2^{\ell-1}}\nabla g_{\eta^{(m)}}(x,\xi),\\
& \overline{g_{\bullet}(x,\xi)}^{\ell-1,(b)} = \frac{1}{2^{\ell-1}}\sum_{m=2^{\ell-1}+1}^{2^{\ell}}g_{\eta^{(m)}}(x,\xi),\quad \overline{\nabla g_{\bullet}(x,\xi)}^{\ell-1,(b)} = \frac{1}{2^{\ell-1}}\sum_{m=2^{\ell-1}+1}^{2^{\ell}}\nabla g_{\eta^{(m)}}(x,\xi),\\
& \psi_{\ell-1}^{(a)}(x) =\left(\overline{\nabla g_{\bullet}(x,\xi)}^{\ell-1,(a)}\right)^{\top}\nabla f_{\xi}\left(\overline{g_{\bullet}(x,\xi)}^{\ell-1,(a)}\right), \quad \mathrm{and}\\
& \psi_{\ell-1}^{(b)}(x) =\left(\overline{\nabla g_{\bullet}(x,\xi)}^{\ell-1,(b)}\right)^{\top}\nabla f_{\xi}\left(\overline{g_{\bullet}(x,\xi)}^{\ell-1,(b)}\right) .
\end{align*}
Here it is important to note that $\psi_{\ell-1}^{(a)}(x)$ and $\psi_{\ell-1}^{(b)}(x)$ are conditionally independent and can be regarded as two different realizations of $\psi_{\ell-1}(x)$ for a common $\xi$, but both are highly correlated with $\psi_{\ell}(x)$ because of the shared use of inner variables. Now let us define
\[ \Delta \psi_{\ell}(x):=\begin{cases} \psi_0(x) & \text{if $\ell=0$,} \\ \psi_{\ell}(x)-\frac{1}{2}\left( \psi_{\ell-1}^{(a)}(x)+\psi_{\ell-1}^{(b)}(x)\right) & \text{otherwise.}\end{cases} \]
Despite of the correlation mentioned above, it follows from the linearity of expectation that
\[ \EE\left[\Delta \psi_{\ell}(x)\right]=\EE\left[\psi_{\ell}(x)\right]-\frac{1}{2}\left( \EE\left[\psi_{\ell-1}^{(a)}(x)\right]+\EE\left[\psi_{\ell-1}^{(b)}(x)\right]\right)=\EE\left[\psi_{\ell}(x)-\psi_{\ell-1}(x)\right], \]
for $\ell>0$, so that we have
\begin{align}\label{eq:telescope}
    \nabla F(x)=\lim_{L\to \infty}\sum_{\ell=0}^{L}\EE\left[\Delta \psi_{\ell}(x)\right].
\end{align}

Before moving on to construction of our MLMC gradient estimator, let us mention that the idea of coupling the successive approximation levels $\psi_{\ell}$ and $\psi_{\ell-1}$ in this way goes back to the work of Giles and Szpruch \cite{GS14} and then has been explored intensively to estimate nested expectations efficiently \cite{BHR15,GG19,GHI20,HGGT20}. The key property that $\Delta \psi_{\ell}(x)$ satisfies is given by the following \emph{antithetic} equalities:
\begin{align}
 \overline{g_{\bullet}(x,\xi)}^{\ell} & = \frac{1}{2}\left( \overline{g_{\bullet}(x,\xi)}^{\ell-1,(a)}+\overline{g_{\bullet}(x,\xi)}^{\ell-1,(b)}\right), \label{eq:antithetic1}\\
 \overline{\nabla g_{\bullet}(x,\xi)}^{\ell} & = \frac{1}{2}\left( \overline{\nabla g_{\bullet}(x,\xi)}^{\ell-1,(a)}+\overline{\nabla g_{\bullet}(x,\xi)}^{\ell-1,(b)}\right).\label{eq:antithetic2}
\end{align}
By exploiting these equalities, it enables us to show that our MLMC gradient estimator is unbiased and also can have finite variance and finite expected cost under some conditions.

\subsection{MLMC estimators}
In constructing the original MLMC estimator from \cite{G08}, we first approximate the infinite series \eqref{eq:telescope} by a partial sum and then estimate each of the terms independently. Given $L\in \NN\cup \{0\}$, by using $N_{\ell}$ i.i.d.\ copies of $\psi_{\ell}(x)$ for $0\leq \ell\leq L$, denoted by $\psi_{\ell}^{(1)}(x),\ldots,\psi_{\ell}^{(N_{\ell})}(x)$, the MLMC estimator is given by
\[ \sum_{\ell=0}^{L}\frac{1}{N_{\ell}}\sum_{m=1}^{N_{\ell}}\Delta \psi_{\ell}^{(m)}(x).\]
Even with a large $L$, the MLMC estimator is generally biased as we have
\begin{align*}
    \EE\left[ \sum_{\ell=0}^{L}\frac{1}{N_{\ell}}\sum_{m=1}^{N_{\ell}}\Delta \psi_{\ell}^{(m)}(x)\right] & =  \sum_{\ell=0}^{L}\frac{1}{N_{\ell}}\sum_{m=1}^{N_{\ell}}\EE\left[\Delta \psi_{\ell}^{(m)}(x)\right] \\
    & =\sum_{\ell=0}^{L}\EE\left[\Delta \psi_{\ell}(x)\right]=\EE\left[ \psi_L(x)\right]\neq \lim_{\ell\to \infty}\EE\left[ \psi_{\ell}(x)\right] = \nabla F(x).
\end{align*} 
The mean squared error can be decomposed into the variance and the squared bias:
\begin{align*}
    & \EE\left[ \left\| \sum_{\ell=0}^{L}\frac{1}{N_{\ell}}\sum_{m=1}^{N_{\ell}}\Delta \psi_{\ell}^{(m)}(x)-\nabla F(x)\right\|_2^2\right] \\
    & = \EE\left[ \left\| \sum_{\ell=0}^{L}\frac{1}{N_{\ell}}\sum_{m=1}^{N_{\ell}}\Delta \psi_{\ell}^{(m)}(x)-\sum_{\ell=0}^{L}\EE\left[\Delta \psi_{\ell}(x)\right]+\EE\left[ \psi_L(x)\right]-\nabla F(x)\right\|_2^2\right] \\
    & = \EE\left[ \left\| \sum_{\ell=0}^{L}\left(\frac{1}{N_{\ell}}\sum_{m=1}^{N_{\ell}}\Delta \psi_{\ell}^{(m)}(x)-\EE\left[\Delta \psi_{\ell}(x)\right]\right)\right\|_2^2\right]+\left\| \EE\left[ \psi_L(x)\right]-\nabla F(x)\right\|_2^2\\
    & = \sum_{\ell=0}^{L}\frac{\EE\left[\left\|\Delta \psi_{\ell}(x)-\EE\left[\Delta \psi_{\ell}(x)\right]\right\|_2^2\right]}{N_{\ell}}+\left\| \EE\left[ \psi_L(x)\right]-\nabla F(x)\right\|_2^2,
\end{align*} 
where we denote the 2-norm of a vector by $\| \cdot\|_2$. Note that the last equality follows from the fact that mutually independent variables are considered between different levels $0\leq \ell\leq L$. The result originally proven in \cite[Theorem~3.1]{G08} can be adapted to our current context straightforwardly as follows, which shows that the MLMC estimator is quite efficient in estimating the pointwise gradient.
\begin{theorem}\label{thm:basic_MLMC}
Assume that there exist constants $\alpha,\beta,c_1,c_2>0$ such that $\alpha\geq \min(\beta,1)/2$,
\begin{align}\label{eq:MLMC_assum}
    \left\| \EE\left[ \psi_{\ell}(x)\right]-\nabla F(x)\right\|_2\leq c_12^{-\alpha \ell} \quad \mathrm{and}\quad \EE\left[\left\|\Delta \psi_{\ell}(x)\right\|_2^2\right]\leq c_22^{-\beta \ell},
\end{align}
for all $\ell\in \NN$. Then, for any given accuracy $\varepsilon<e^{-1}$, there exist $L$ and $N_0,N_1,\ldots,N_L$ such that the MLMC estimator has a mean squared error bound
\[ \EE\left[ \left\| \sum_{\ell=0}^{L}\frac{1}{N_{\ell}}\sum_{m=1}^{N_{\ell}}\Delta \psi_{\ell}^{(m)}(x)-\nabla F(x)\right\|_2^2\right]\leq \varepsilon^2,\]
with the total computational cost 
\[ \sum_{\ell=0}^{L}N_{\ell}2^{\ell} \leq \begin{cases} c_4\varepsilon^{-2}, & \text{if $\beta>1$,} \\ c_4\varepsilon^{-2}|\log \varepsilon^{-1}|^2, & \text{if $\beta=1$,} \\ c_4\varepsilon^{-2-(1-\beta)/\alpha}, & \text{if $\beta<1$,} \end{cases} \]
for a constant $c_4>0$ independent of $\varepsilon$.
\end{theorem}

As discussed, for instance, in \cite[Section~2.1]{GG19}, in any regime whether $\beta>1$, $\beta=1$, or $\beta<1$, the total computational cost of the nested Monte Carlo estimator \eqref{eq:NMC}, which is given by the product $MN$, is only bounded above by $c_5\varepsilon^{-2-1/\alpha}$ for a constant $c_5>0$ to achieve a given accuracy $\varepsilon$. Therefore, the MLMC estimator always has an asymptotically better complexity bound.

In the regime $\beta>1$, we can construct an unbiased gradient estimator with finite variance and expected cost per one outer sample of $\xi$ by applying a randomization technique from \cite{M11,RG15}. Under the assumption \eqref{eq:MLMC_assum}, for any sequence $\omega_0,\omega_1,\ldots>0$ such that $\omega_0+\omega_1+\cdots=1$, it follows from \eqref{eq:telescope} that we have an absolutely convergent series
\begin{align*}
    \nabla F(x)=\sum_{\ell=0}^{\infty}\EE\left[\Delta \psi_{\ell}(x)\right]=\sum_{\ell=0}^{\infty}\frac{\EE\left[\Delta \psi_{\ell}(x)\right]}{\omega_\ell}\omega_\ell.
\end{align*}
This means that, for a level $\ell\geq 0$ randomly chosen with probability $\omega_{\ell}$, 
\[ \frac{\Delta \psi_{\ell}(x)}{\omega_{\ell}}\]
is an unbiased gradient estimator. More generally, for any $N\in \NN$, with $\ell^{(1)},\ldots,\ell^{(N)}\geq 0$ being independently and randomly chosen from the discrete distribution defined by the sequence $\omega_0,\omega_1,\ldots$ and $\Delta \psi_{\ell^{(1)}},\ldots,\Delta \psi_{\ell^{(N)}}$ being mutually independent, we can construct an unbiased estimator by
\begin{align}\label{eq:unbiased_MLMC}
    \frac{1}{N}\sum_{n=1}^{N}\frac{\Delta \psi_{\ell^{(n)}}(x)}{\omega_{\ell^{(n)}}}.
\end{align}
In the literature, an estimator of this form is called a \emph{single term estimator} \cite{RG15,V18}. Although we can also define the corresponding \emph{independent/couple sum estimators}, we focus on the above estimator for the sake of simplicity. Moreover, one can use the stratified sampling or some other techniques for generating $\ell^{(1)},\ldots,\ell^{(N)}\geq 0$, which still lead to unbiased MLMC estimators \cite{V18}.

Here, the variance per one outer sample is given by
\[    \EE\left[ \left\|  \frac{\Delta \psi_{\ell}(x)}{\omega_{\ell}}-\nabla F(x)\right\|_2^2\right] = \sum_{\ell=0}^{\infty}  \frac{\EE\left[ \left\|\Delta \psi_{\ell}(x)\right\|_2^2\mid \ell \right]}{\omega_{\ell}}-\|\nabla F(x)\|_2^2 \leq c_2\sum_{\ell=0}^{\infty}  \frac{2^{-\beta \ell}}{\omega_{\ell}}, \]
while the expected computational cost is 
\[ \sum_{\ell=0}^{\infty}\omega_{\ell}2^{\ell}.\]
Both are finite by choosing \tg{$\omega_{\ell}\propto 2^{-\tau \ell}$ for any $\tau\in (1,\beta)$} if $\beta>1$. This way, the CSO problems \eqref{eq:CSO} can be solved, for instance, by the first-order stochastic gradient descent
\begin{align}\label{eq:SGD}
    x_{t+1}=x_t-\frac{\gamma_t}{N}\sum_{n=1}^{N}\frac{\Delta \psi_{\ell^{(n)}}(x_t)}{\omega_{\ell^{(n)}}},\quad t=0,1,\ldots,
\end{align}
with a initial point $x_0\in \Xcal$ and a sequence of step-sizes $\gamma_0,\gamma_1,\ldots>0$, due to Robbins and Monro \cite{RM51}. There are many variants of this Robbins-Monro algorithm, such as Polyak-Ruppert averaging \cite{P90,R91}, the stochastic counterpart of Nesterov's acceleration \cite{N83}, and adaptation methods \cite{DHS11,KB15,RKK19}, in which our unbiased MLMC gradient estimator can be used.

\begin{remark}
Since our MLMC gradient estimator is unbiased and has finite variance, the standard theory from stochastic optimization directly applies to the CSO problems. For example, if the domain $\Xcal$ is convex, $F$ is strongly convex and differentiable, and $\EE\left[ \left\| \Delta \psi_{\ell}(x)/\omega_{\ell}\right\|_2^2\right]<\infty$ for any $x\in \Xcal$, our estimate $x_t$ obtained by the recursion \eqref{eq:SGD} converges to the optimal $x^*$ with the mean squared error of $O(t^{-1})$ as long as the step-sizes satisfy
\[ \sum_{t=0}^{\infty}\gamma_t=\infty\quad \text{and}\quad \sum_{t=0}^{\infty}\gamma^2_t<\infty. \]
We refer to \cite{KY03,SDP09,L20} for more information and also \cite{GLQSSR19} and the references cited therein for recent advances in this research direction. Moreover, a sequence of the works by Hu et al.\ \cite{HCH20,HZCH20,HCH21} studies the computational complexity of both the standard nested Monte Carlo and MLMC gradient estimators for the CSO problems carefully, which concludes a superiority of the MLMC estimator (with $\beta=1$) for several convexity classes. Our theoretical result below does not improve the complexity for the CSO problems but enables to avoid a delicate bias-variance trade-off in computation without any need to specify the maximum level $L$ (and $N_0,\ldots,N_L$) for the original MLMC estimator.
\end{remark}

In the light of the above remark, it is important to show when the condition $\beta>1$ holds, which we work on in the next subsection. 

\subsection{Theoretical analysis}
As a main result of this paper, we prove a bound on the variance for our unbiased MLMC estimator \tg{in this subsection}. In what follows, we denote the $p$-norm of a vector $v$ by $\|v\|_p$. For a matrix $A$, we use the same notation to represent the $L_{p,p}$-norm, i.e., $\|A\|_p=\left(\sum_i\sum_ja_{ij}^p\right)^{1/p}$. Thus, $\|A\|_2$ is nothing but the Frobenius norm of $A$. 

\begin{theorem}\label{thm:main1}
Assume that $\nabla f_{\xi}$ is Lipschitz continuous and $\nabla^{\top}\nabla f_{\xi}$ is $\rho$-H\"{o}lder continuous, i.e., assume that there exist  $\lambda_1,\lambda_2>0$ and $0<\rho\leq 1$ such that
\begin{align}\label{assump1} 
\left\|\nabla f_{\xi}(g_{\eta_1})-\nabla f_{\xi}(g_{\eta_2})\right\|_2 \leq \lambda_1\|g_{\eta_1}-g_{\eta_2}\|_2,
\end{align}
and 
\begin{align}\label{assump2} 
 \left\|\nabla^{\top}\nabla f_{\xi}(g_{\eta_1})-\nabla^{\top}\nabla f_{\xi}(g_{\eta_2})\right\|_2 \leq  \lambda_2\|g_{\eta_1}-g_{\eta_2}\|_2^{\rho},
\end{align}
hold for any $g_{\eta_1}$ and $g_{\eta_2}$. Additionally assume that there exist $p\geq 4$ and $q\geq 4(1+\rho)$ such that 
\[ \EE_{\xi}\EE_{\eta\mid \xi}\left[\left\| g_{\eta}\right\|_{p}^{p}\right]<\infty \quad \text{and}\quad  \EE_{\xi}\EE_{\eta\mid \xi}\left[\left\|\nabla g_{\eta}\right\|_{q}^{q}\right]<\infty. \]
Then \tg{we have
$\EE\left[\left\|\Delta \psi_{\ell}(x)\right\|_2^2\right]\leq c_2 2^{-(1+\rho)\ell}$ for all $\ell\in \NN$, where $c_2>0$ is a constant independent of $\ell$.}
\end{theorem}

\tg{We shall give a proof of the theorem in \ref{app:proof_main1}. In reference to Theorem~\ref{thm:basic_MLMC}, the result of Theorem~\ref{thm:main1} directly means that $\beta=1+\rho>1$, so that we can construct an unbiased MLMC gradient estimator with finite variance and expected cost as follows.}

\tg{\begin{corollary}\label{cor:main1}
Under the same assumptions as in Theorem~\ref{thm:main1}, by choosing $\omega_{\ell}\propto 2^{-\tau \ell}$ for $\tau\in (1,1+\rho)$, the gradient estimator \eqref{eq:unbiased_MLMC} is unbiased and has finite variance and expected cost for any $N\in \NN$.
\end{corollary}}

\section{Another gradient estimation for a special case}\label{sec:special}
Finally, as one of the important class of CSO problems, let us consider the situation where $u(\xi)\in \RR$ represents an observable variable, given as a function of the hidden variable $\xi$, and $f_{\xi}$ measures the squared error between $u(\xi)$ and our predictor $\EE_{\eta\mid \xi}\left[ g_{\eta}\left(x, \xi\right)\right]$ for a function $g_{\eta}(\cdot,\xi): \RR^d\to \RR$. That is, the objective function $F(x)$ is given by
\[ F(x)=\EE_{\xi}\left[ \left(u(\xi)- \EE_{\eta\mid \xi}\left[ g_{\eta}\left(x, \xi\right)\right]\right)^2\right]. \]
The CSO problems of this form appear, for instance, in instrumental variable regression \cite{HLLT17,SSG19}, \tg{as described in Example~\ref{exm:iv}.} In a slightly different context, estimating nested expectations of this form has been investigated in \cite{G17,WLZ17} and also mentioned as a special case in \cite{RCYWW18}. However, estimating its gradient $\nabla F(x)$ has not been studied well in the literature.

\subsection{Unbiased gradient estimators}
Let $\eta'$ be the i.i.d.\ copy of $\eta$ conditioned on $\xi$. The gradient can be calculated as
\begin{align*}
    \nabla F(x) & =  -2\, \EE_{\xi}\left[  \left(u(\xi)- \EE_{\eta\mid \xi}\left[ g_{\eta}\left(x, \xi\right)\right]\right)\EE_{\eta\mid \xi}\left[ \nabla g_{\eta}\left(x, \xi\right)\right]\right] \\
    & = -2\, \EE_{\xi}\EE_{\eta,\eta'\mid \xi}\left[  \left(u(\xi)- g_{\eta}\left(x, \xi\right)\right) \nabla g_{\eta'}\left(x, \xi\right)\right].
\end{align*}
With these representations, it is possible to construct several unbiased gradient estimators as follows. 

For positive integers $M,N$, let $\xi^{(1)},\ldots,\xi^{(N)}$ be the i.i.d.\ copies of $\xi$, and for each $n=1,\ldots,N$, $\eta^{(1,n)},\ldots,\eta^{(M,n)}$ and $\eta'^{(1,n)},\ldots,\eta'^{(M,n)}$ be the i.i.d.\ copies of $\eta$ conditioned on $\xi=\xi^{(n)}$. Then, because of the conditional independence between $\eta^{(1,n)},\ldots,\eta^{(M,n)}$ and $\eta'^{(1,n)},\ldots,\eta'^{(M,n)}$, the gradient $\nabla F(x)$ can be estimated without any bias by
\begin{align*}
    \widehat{\nabla F}^{(1)}_{M,N}(x) & = -\frac{2}{N}\sum_{n=1}^{N} \left(u(\xi^{(n)})-\frac{1}{M}\sum_{m=1}^{M}g_{\eta^{(m,n)}}(x, \xi^{(n)})\right)\frac{1}{M}\sum_{m=1}^{M}\nabla g_{\eta'^{(m,n)}}(x, \xi^{(n)}),\\
    \widehat{\nabla F}^{(2)}_{M,N}(x) & = -\frac{1}{N}\sum_{n=1}^{N} \left[\left(u(\xi^{(n)})-\frac{1}{M}\sum_{m=1}^{M}g_{\eta^{(m,n)}}(x, \xi^{(n)})\right)\frac{1}{M}\sum_{m=1}^{M}\nabla g_{\eta'^{(m,n)}}(x, \xi^{(n)})\right.\\
    & \quad \quad \quad \quad \quad \quad \left. + \left(u(\xi^{(n)})-\frac{1}{M}\sum_{m=1}^{M}g_{\eta'^{(m,n)}}(x, \xi^{(n)})\right)\frac{1}{M}\sum_{m=1}^{M}\nabla g_{\eta^{(m,n)}}(x, \xi^{(n)})\right].
\end{align*}
We note that the first estimator $\widehat{\nabla F}^{(1)}_{M,N}(x)$ has been used in \cite{HLLT17}. 
\tg{It is worth observing here that, apart from the factor $1/2$, the two halves of $\widehat{\nabla F}^{(2)}_{M,N}(x)$ are both identical in distribution to $\widehat{\nabla F}^{(1)}_{M,N}(x)$. Thus, denoting the correlation of the $i$-th component of the two halves of $\widehat{\nabla F}^{(2)}_{M,N}(x)$ by $\rho_i\in [-1,1]$, we have
\begin{align*}
    \EE\left[ \left\|\widehat{\nabla F}^{(2)}_{M,N}(x)-\nabla F(x)\right\|_2^2\right] & \leq \frac{1}{2}\left( 1+\max_{i}\rho_i\right)\EE\left[ \left\|\widehat{\nabla F}^{(1)}_{M,1}(x)-\nabla F(x)\right\|_2^2\right] \\
    & \leq \EE\left[ \left\|\widehat{\nabla F}^{(1)}_{M,1}(x)-\nabla F(x)\right\|_2^2\right].
\end{align*}
This means that the variance of $\widehat{\nabla F}^{(1)}_{M,N}(x)$ is always larger than or equal to that of $\widehat{\nabla F}^{(2)}_{M,N}(x)$. In the proof of Theorem~\ref{thm:main2} below, we evaluate the difference between these variances more precisely.}

Here we provide yet another unbiased estimator for this special case. For a start, \tg{let us consider a problem of estimating $(\EE[X])^2$ for a real-valued random variable $X$ with variance $\VV[X]<\infty$. Denoting the $M$ i.i.d.\ copies of $X$ by $X^{(1)},\ldots,X^{(M)}$, an estimator
\[ \left(\frac{1}{M}\sum_{m=1}^{M}X^{(m)}\right)^2 \]
is biased, as we have
\[ \EE\left[ \left(\frac{1}{M}\sum_{m=1}^{M}X^{(m)}\right)^2\right] = \left( \EE\left[ \frac{1}{M}\sum_{m=1}^{M}X^{(m)}\right] \right)^2 + \VV\left[ \frac{1}{M}\sum_{m=1}^{M}X^{(m)}\right] = \left( \EE\left[X\right] \right)^2 + \frac{\VV\left[ X\right]}{M}.
\]
It is well known, however, that an unbiased estimator for the variance $\VV[X]$ is given by
\[ \frac{1}{M-1}\sum_{m=1}^{M}\left( X^{(m)}-\frac{1}{M}\sum_{m'=1}^{M}X^{(m')}\right)^2, \]
so that we can construct a bias-corrected estimator for $(\EE[X])^2$ as
\[ \left(\frac{1}{M}\sum_{m=1}^{M}X^{(m)}\right)^2-\frac{1}{M(M-1)}\sum_{m=1}^{M}\left( X^{(m)}-\frac{1}{M}\sum_{m'=1}^{M}X^{(m')}\right)^2.\]}

\tg{Now, to estimate $F(x)$, define
\[ \widehat{F}_{M,N}(x) = \frac{1}{N}\sum_{n=1}^{N}\left( u(\xi^{(n)})-\frac{1}{M}\sum_{m=1}^{M}g_{\eta^{(m,n)}}(x, \xi^{(n)})\right)^2.\]
This is a biased estimator for the objective function $F(x)$, where the bias comes from estimating the term $\EE_{\xi}\left[ \left(\EE_{\eta\mid \xi}\left[ g_{\eta}\left(x, \xi\right)\right]\right)^2\right]$. Following the above argument, we can give its bias-corrected version
\begin{align}\label{eq:special_unbiased}
    \widehat{F}^*_{M,N}(x) = \widehat{F}_{M,N}(x)-\frac{1}{N}\sum_{n=1}^{N}\frac{1}{M(M-1)}\sum_{m=1}^{M}\left( g_{\eta^{(m,n)}}(x, \xi^{(n)})-\frac{1}{M}\sum_{m'=1}^{M}g_{\eta^{(m',n)}}(x, \xi^{(n)})\right)^2,
\end{align}
for any $M\geq 2$. That is, for any positive integers $M\geq 2$ and $N$, we have $\EE[\widehat{F}^*_{M,N}(x)]=F(x)$.}

Motivated by $\widehat{F}^*_{M,N}(x)$, we introduce the third unbiased gradient estimator 
\begin{align*}
    \widehat{\nabla F}^{(3)}_{M,N}(x) & = \nabla \widehat{F}^*_{M,N}(x) \\
    & = \frac{1}{N}\sum_{n=1}^{N}\left[-2\left( u(\xi^{(n)})-\frac{1}{M}\sum_{m=1}^{M}g_{\eta^{(m,n)}}(x, \xi^{(n)})\right)\frac{1}{M}\sum_{m=1}^{M}\nabla g_{\eta^{(m,n)}}(x, \xi^{(n)})\right. \\
    & \quad \quad \quad \quad \quad \quad \left.-\frac{2}{M(M-1)}\sum_{m=1}^{M}\left( g_{\eta^{(m,n)}}(x, \xi^{(n)})-\frac{1}{M}\sum_{m'=1}^{M}g_{\eta^{(m',n)}}(x, \xi^{(n)})\right)\right.\\
    & \quad \quad \quad \quad \quad \quad \quad \quad \quad \quad \quad \quad \times \left. \left( \nabla g_{\eta^{(m,n)}}(x, \xi^{(n)})-\frac{1}{M}\sum_{m'=1}^{M}\nabla g_{\eta^{(m',n)}}(x, \xi^{(n)})\right)\right]. 
\end{align*}
Here the unbiasedness of $\widehat{\nabla F}^{(3)}_{M,N}(x)$ can be proven \tg{by elementary calculations,} so we omit the proof. Note that $\widehat{\nabla F}^{(3)}_{M,N}(x)$ with $M=2$ coincides with the the second estimator $\widehat{\nabla F}^{(2)}_{M,N}(x)$ with $M=1$.

\begin{remark}
The first two estimators $\widehat{\nabla F}^{(1)}_{M,N}(x)$ and $\widehat{\nabla F}^{(2)}_{M,N}(x)$ require $2M$ inner conditional samples of $\eta$ for each outer sample of $\xi$, whereas the last estimator $\widehat{\nabla F}^{(3)}_{M,N}(x)$ does $M$ inner samples of $\eta$. Thus, from the point of view of the computational cost, it is fair to compare the variance of $\widehat{\nabla F}^{(1)}_{M,N}(x)$,  $\widehat{\nabla F}^{(2)}_{M,N}(x)$ and $\widehat{\nabla F}^{(3)}_{2M,N}(x)$ for given $M$.
\end{remark}

\subsection{Theoretical analysis}
Our theoretical result for this special case is given as follows, which shows that the third estimator has the minimum variance among the three under a fixed cost.

\begin{theorem}\label{thm:main2}
    Assume that $\EE_{\xi}\left[ (u(\xi))^2\right],\EE_{\xi}\EE_{\eta\mid \xi}\left[\left\| g_{\eta}\right\|_2^2\right],\EE_{\xi}\EE_{\eta\mid \xi}\left[\left\|\nabla g_{\eta}\right\|_2^2\right]<\infty$. \tg{For any positive integers $M$ and $N$, we have
    \[ \EE\left[ \left\|\widehat{\nabla F}^{(1)}_{M,N}(x)-\nabla F(x)\right\|_2^2\right] \geq \EE\left[ \left\|\widehat{\nabla F}^{(2)}_{M,N}(x)-\nabla F(x)\right\|_2^2\right]\geq \EE\left[ \left\|\widehat{\nabla F}^{(3)}_{2M,N}(x)-\nabla F(x)\right\|_2^2\right]. \]}
\end{theorem}
\noindent \tg{A proof of the theorem is given in \ref{app:proof_main2}.}

\section{Numerical experiments}\label{sec:numerical}
\subsection{Invariant logistic regression}
In order to confirm the effectiveness of our unbiased MLMC gradient estimator, we \tg{first} examine it with a simple test case on invariant logistic regression studied in \cite{HCH20,HZCH20}. The objective function is given by
\[ F(x)=\EE_{\xi=(a,b)}\left[ \log\left(1+\exp\left(-b\EE_{\eta\mid \xi}\left[ \eta\right]^{\top}x\right)\right)\right],\]
where $\xi$ consists of $a\in \RR^d$ and $b\in \{\pm 1\}$ with $a$ being the random feature vector and $b$ the corresponding label, and $\eta$ is a randomly perturbed observation of the vector $a$. In what follows, we consider the case where $d=10$, $a\sim \Ncal(\bszero,\sigma_{\xi}^2I_{d})$ with $\bszero$ denoting the vector of all zeros and $I_d$ the $d\times d$ identity matrix, and $\eta\mid \xi \sim \Ncal(a,\sigma_{\eta}^2I_{d})$. The label $b$ is given by
\[ b=\begin{cases} 1 & \text{if $a^{\top}x^*>0,$}\\ -1 & \text{otherwise,}\end{cases} \]
with the fixed vector $x^*=(1,2,\ldots,d)$. For simplicity, we set $\sigma_{\xi}^2=\sigma_{\eta}^2=1$.

\begin{figure}
    \centering
    \includegraphics[width=0.4\textwidth]{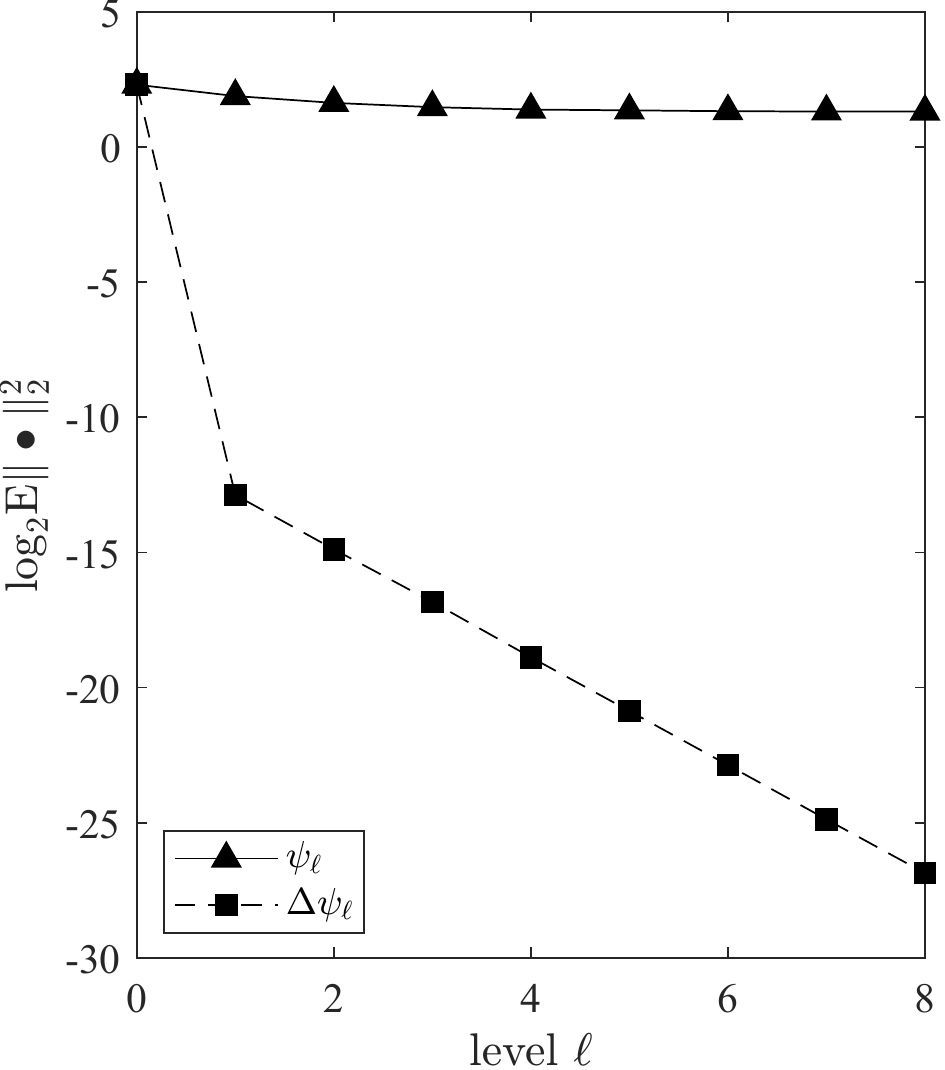}
    \includegraphics[width=0.4\textwidth]{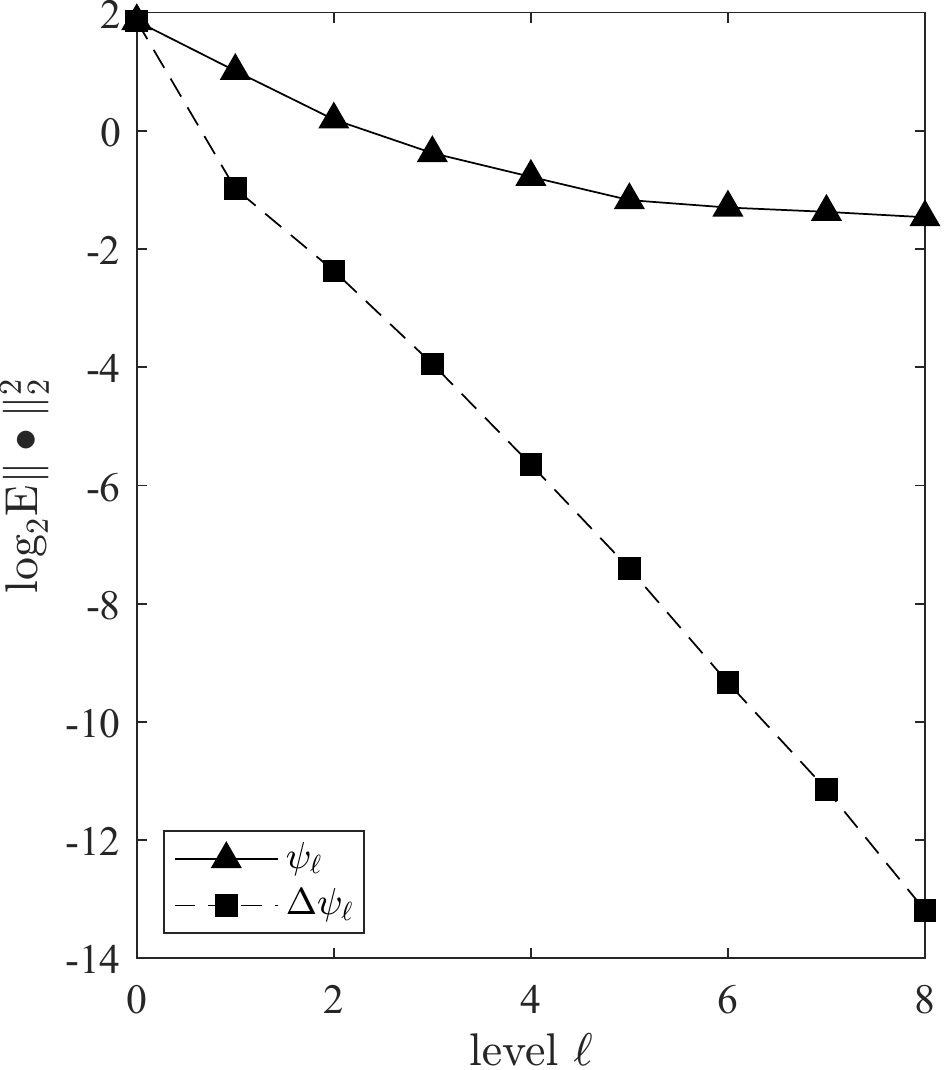}
    \caption{The expected squared 2-norms of the variables $\psi_{\ell}(x)$ and $\Delta \psi_{\ell}(x)$ for the logistic regression model at a randomly generated point from $\Ncal(0,10^{-4}I_d)$ (left) and at a point found after the Robbins-Monro iterations (right).}
    \label{fig:mlmc_coupling}
\end{figure}

The left panel of Figure~\ref{fig:mlmc_coupling} shows the convergence behavior of $\EE\left[\left\|\Delta \psi_{\ell}(x)\right\|_2^2\right]$ as a function of $\ell$ for one randomly generated point $x$ from $\Ncal(0,10^{-4}I_d)$. As comparison, the corresponding result for $\EE\left[\left\|\psi_{\ell}(x)\right\|_2^2\right]$ is also shown in the panel. As it is not possible to compute the exact values of $\EE\left[\left\|\Delta \psi_{\ell}(x)\right\|_2^2\right]$ and $\EE\left[\left\|\psi_{\ell}(x)\right\|_2^2\right]$, we estimate them by the Monte Carlo averages over $10^4$ i.i.d.\ samples of $\Delta \psi_{\ell}(x)$ and $\psi_{\ell}(x)$, respectively, for each $\ell=0,1,\ldots,8$. We can see that $\EE\left[\left\|\psi_{\ell}(x)\right\|_2^2\right]$ takes an almost constant value for the considered range of $\ell$, whereas $\EE\left[\left\|\Delta \psi_{\ell}(x)\right\|_2^2\right]$ decays geometrically fast. The linear regression for $\log_2\left(\EE\left[\left\|\psi_{\ell}(x)\right\|_2^2\right]\right)$ with $1\leq \ell\leq 8$ provides an estimation of $\beta$ as 2.00, which agrees well with our theoretical result.

According to \eqref{eq:SGD}, we minimize $F(x)$ by the \tg{Robbins}-Monro algorithm with the unbiased MLMC gradient estimator. We set $N=1$ and $\omega_{\ell}=(1-2^{-3/2})2^{-3\ell/2}$ for all $\ell\in \NN\cup\{0\}$. This means that the expected cost of our MLMC gradient estimator per one iteration step is given by
\[ \sum_{\ell=0}^{\infty}\omega_{\ell}2^{\ell}=\frac{1-2^{-3/2}}{1-2^{-1/2}}\approx 2.21. \]
The initial point $x_0$ is randomly generated from $\Ncal(0,10^{-4}I_d)$, the step-sizes $\gamma_0,\gamma_1,\ldots$ are fixed to be $10^{-4}$, and the \tg{Robbins}-Monro iteration terminates after $T$ steps such that the total expected cost $T(1-2^{-3/2})/(1-2^{-1/2})$ does not exceed $10^6$. We conduct 10 independent runs under these setups. As comparison, we also consider minimizing $F(x)$ by the \tg{Robbins}-Monro algorithm with the biased standard Monte Carlo gradient estimator \eqref{eq:NMC}, where we fix $N=1$ and choose $M=1,2,4,8,16$. For each $M$, the \tg{Robbins}-Monro iteration terminates after $T$ steps such that the total cost $MT$ does not exceed $10^6$.

\begin{figure}
    \centering
    \includegraphics[width=0.8\textwidth]{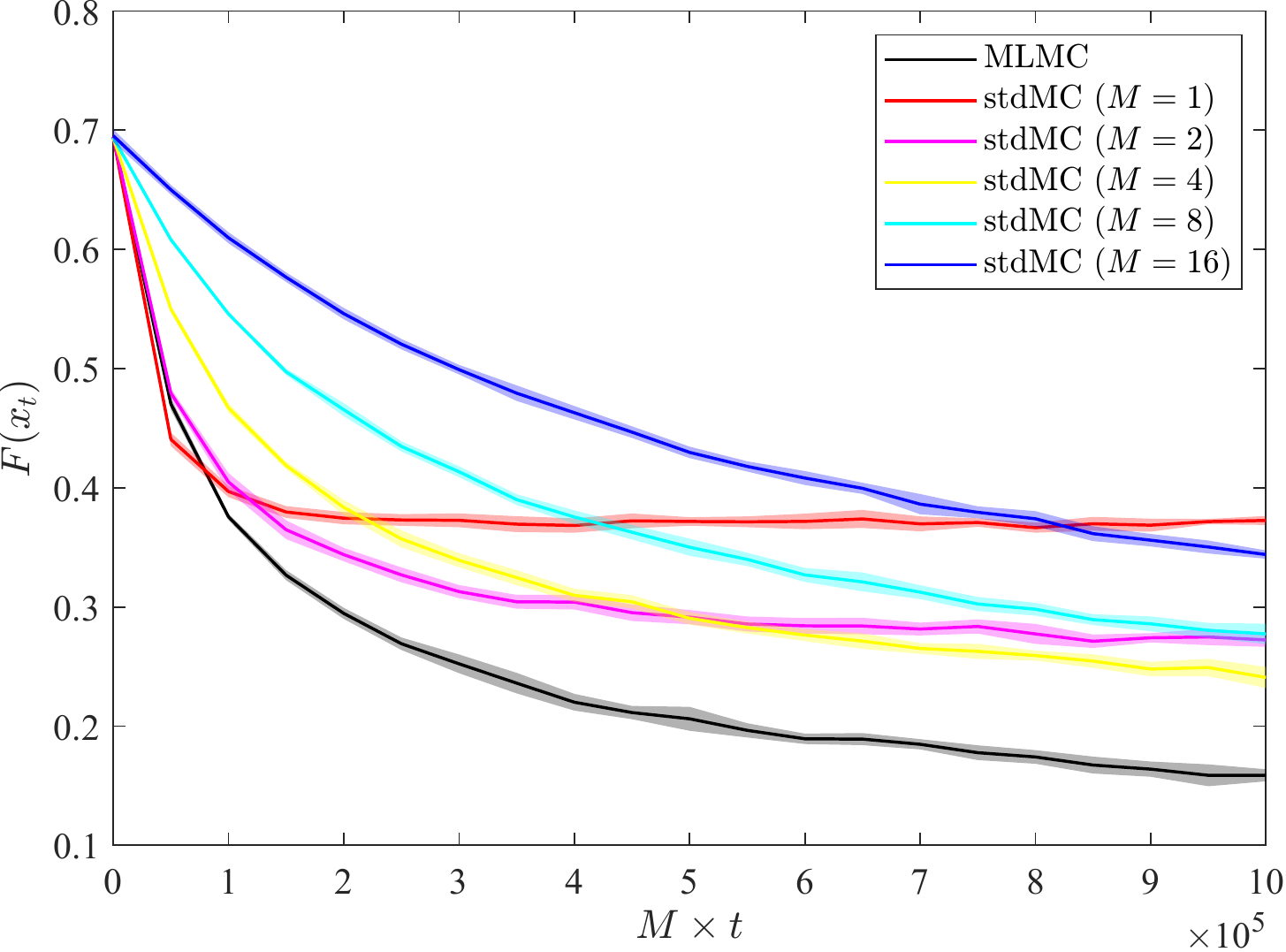}
    \caption{The convergence behavior of the objective function values for the logistic regression model by the \tg{Robbins}-Monro algorithm with several gradient estimators. For each gradient estimator, the solid line and the shaded area represent the average and its standard error estimated from 10 independent runs, respectively. Here, stdMC stands for the standard Monte Carlo gradient estimator.}
    \label{fig:convergence}
\end{figure}

Figure~\ref{fig:convergence} shows the behaviors of the objective function value $F(x_t)$ as a function of the product of the (expected) computational cost per one iteration, $M$, and the number of iteration steps, $t$. For our unbiased MLMC estimator, we simply put $M=(1-2^{-3/2})/(1-2^{-1/2})$. As we cannot compute the exact value of $F(x)$ but the inner conditional expectation $\EE_{\eta\mid \xi}\left[ \eta\right]$ is equal to the outer variable $a$, we estimate $F(x)$ by the simple Monte Carlo average
\[ F(x)\approx \frac{1}{\hat{N}}\sum_{n=1}^{\hat{N}} \log\left(1+\exp\left(-b^{(n)}{a^{(n)}}^{\top}x\right)\right),\]
with $\hat{N}=10^3$, for $(a^{(1)},b^{(1)}),\ldots,(a^{(\hat{N})},b^{(\hat{N})})$ being the i.i.d.\ samples of $\xi$. As 10 independent runs are performed for each gradient estimator, we draw the solid line to represent the average of 10 estimated values and the shaded area is shown to represent the standard error. 

It is interesting to see that the biased standard Monte Carlo gradient estimator with small value of $M$ yields the convergence to a large value of $F$, which clearly shows the disadvantage that the ``bias level" remains the same throughout the iterations. In fact, the biased standard Monte Carlo gradient estimator with $M=1$ for $\nabla F(x)$ coincides the unbiased standard Monte Carlo gradient estimator for
\[ \EE_{\xi}\EE_{\eta\mid \xi}\left[ f_{\xi}\left(  g_{\eta}\left(x, \xi\right)\right)\right], \]
so that the stochastic gradient descent will minimize a wrong objective function. The situation may get improved gradually as $M$ increases, whereas, as can be seen from the figure, the speed of the convergence as a function of $Mt$ becomes slower and we cannot see any significant improvement under our considered setups. In contrast, our MLMC estimator works very well and succeeds in finding much smaller values of $F$.

The right panel of Figure~\ref{fig:mlmc_coupling} shows the convergence behaviors of $\EE\left[\left\|\Delta \psi_{\ell}(x)\right\|_2^2\right]$ and $\EE\left[\left\|\psi_{\ell}(x)\right\|_2^2\right]$ as functions of $\ell$ for the point $x_T$, i.e., the point found after a single run with the unbiased MLMC estimator. Here again, we estimate them by the Monte Carlo averages over $10^4$ i.i.d.\ samples of $\Delta \psi_{\ell}(x)$ and $\psi_{\ell}(x)$, respectively, for each $\ell=0,1,\ldots,8$. We can confirm the geometric decay of $\EE\left[\left\|\Delta \psi_{\ell}(x)\right\|_2^2\right]$, and the linear regression provides an estimation of $\beta$ as 1.75. This implies that, not only at the beginning of the Robbins-Monro iteration, the variance of our MLMC estimator remains bounded with our choice $\omega_{\ell}\propto 2^{-3\ell/2}$, \tg{as we have $1.5\in (1, 1.75)$}.

\subsection{Instrumental variable regression}
\tg{We also examine, for an IV regression problem, the effectiveness of our unbiased MLMC gradient estimator as well as the one with the smallest variance introduced in Section~\ref{sec:special}. We refer to Example~\ref{exm:iv} for the problem setting. Similarly to the exposition in \cite[Section~5.1]{BKS19}, assume that the data is generated by the following process:
\begin{align*}
& Y = f(X)+e+\delta,\qquad X=\frac{Z_1+e}{2}+\gamma,\\
& Z\sim U([-3,3]^2),\qquad e\sim N(0,1),\qquad \gamma,\delta\overset{\text\small\textrm{iid}}{\sim} N(0, 0.1).
\end{align*}
We consider the four cases for the ground truth function $f$ as
\[ f(x)=\sin(x),\quad f(x)=x,\quad f(x)=|x|,\quad f(x)=\bsone_{x\geq 0}.\]
For each case, sampling the outer variable $Y$ is done by following the above generating process with the respective function, whereas $f$ is assumed to be unknown and needs to be optimized by solving \eqref{eq:CSO_IV}. In the usual setting, we only have a finite set of data points for $(X,Y,Z)$, whereas both an infinite supply of the pair $(X,Y,Z)$ and the exact conditional sampling of $X$ given $Z$ are allowed in this experiment. We estimate the unknown $f$ by a neural network consisting of two hidden layers with 50 nodes and the ReLU as an activation function between each layer. This means that, in \eqref{eq:CSO}, $x$ is given by a set of weight parameters of neural network, and we write the corresponding neural network by $g_x$.
}

\tg{By employing the squared loss $L(x,y)=(x-y)^2$, our problem reduces to minimize the objective function
\[ F(x)=\EE_{Y,Z}\left[ \left(Y-\EE_{X\mid Z}\left[ g_x(X)\right]\right)^2\right].\]
We do this by the Robbins-Monro algorithm \eqref{eq:SGD} with the unbiased MLMC gradient estimator. We set $N=1$ and $\omega_{\ell}=(1-2^{-3/2})2^{-3\ell/2}$ for all $\ell\in \NN\cup\{0\}$ again in this example, so that the expected cost per one iteration step is $(1-2^{-3/2})/(1-2^{-1/2})\approx 2.21$. Although we do not describe in detail, we have checked the decay of $\EE\left[\left\|\Delta \psi_{\ell}(x)\right\|_2^2\right]$ for various cases of $f$ and $x$ by following the same way with what is done in Figure~\ref{fig:mlmc_coupling}, and the values of $\beta$ have been consistently estimated to be almost $2.0$. This pre-experiment motivates the above choice for $\omega_{\ell}$. We also minimize $F(x)$ by the Robbins-Monro algorithm in conjunction the third unbiased gradient estimator proposed in Section~\ref{sec:special}. As proven in Theorem~\ref{thm:main2}, $\widehat{\nabla F}^{(3)}_{M,N}$ has the smallest variance among the proposed three estimators under a fixed cost. We fix $N=1$ and choose $M=2,3,4,8$ to make the cost per one iteration step comparable to that of the MLMC estimator.}

\begin{figure}
    \centering
    \begin{minipage}{0.45\textwidth}
    \centering
    \includegraphics[width=\textwidth]{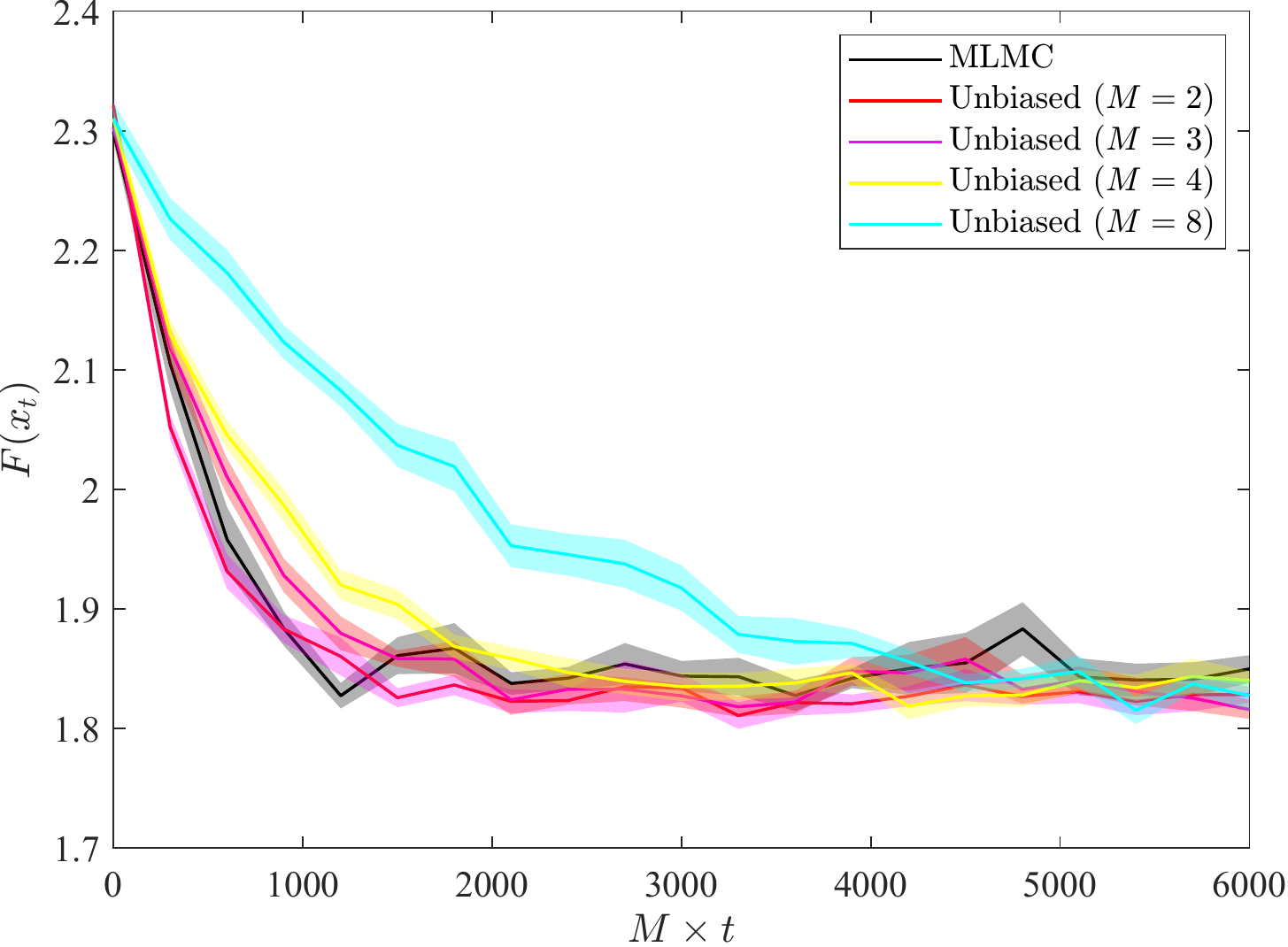}
    \subcaption{$f(x)=\sin(x)$}
    \end{minipage}
    \begin{minipage}{0.45\textwidth}
    \centering
    \includegraphics[width=\textwidth]{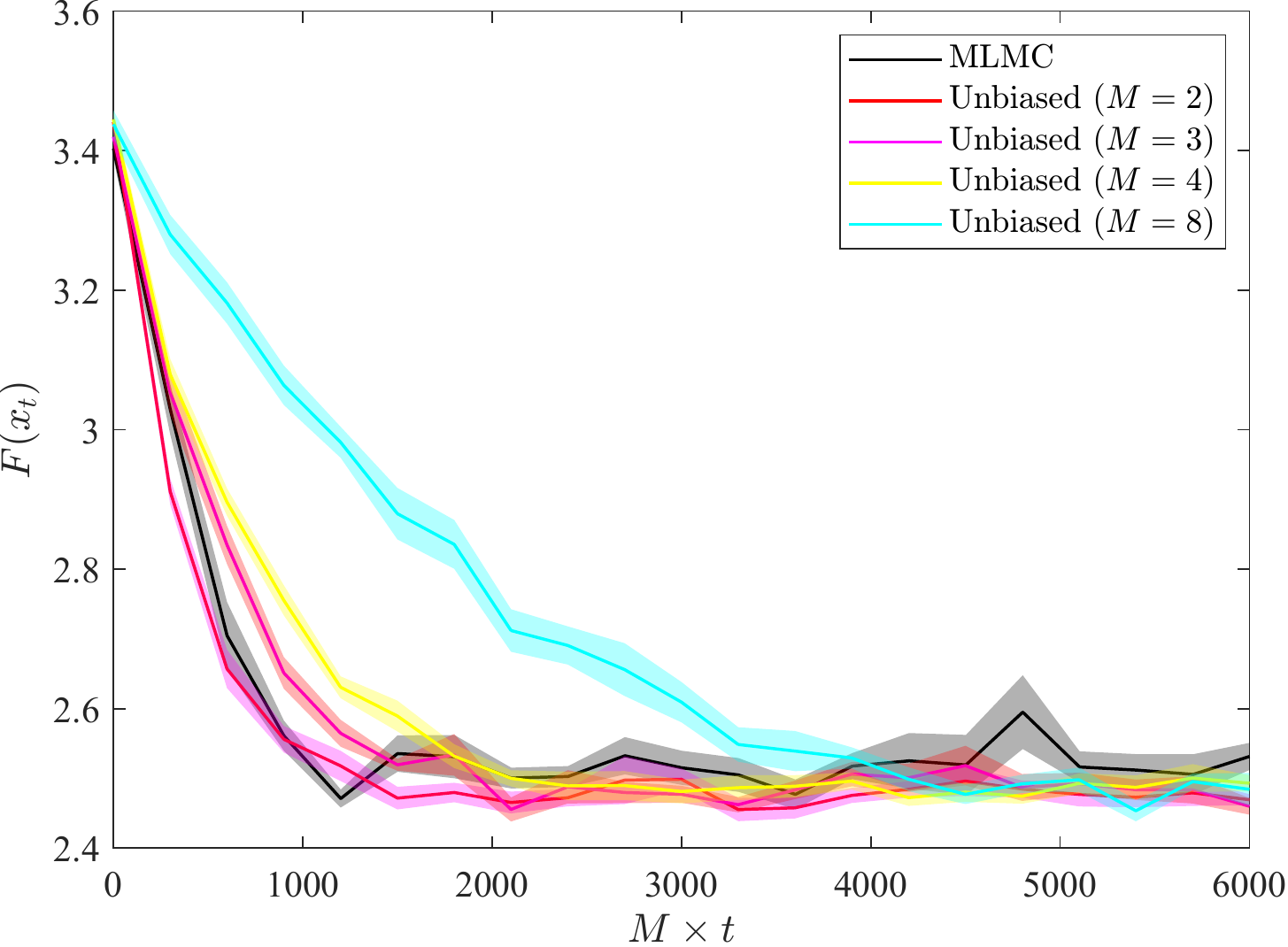}
    \subcaption{$f(x)=x$}
    \end{minipage}\\[5pt]
    \begin{minipage}{0.45\textwidth}
    \centering
    \includegraphics[width=\textwidth]{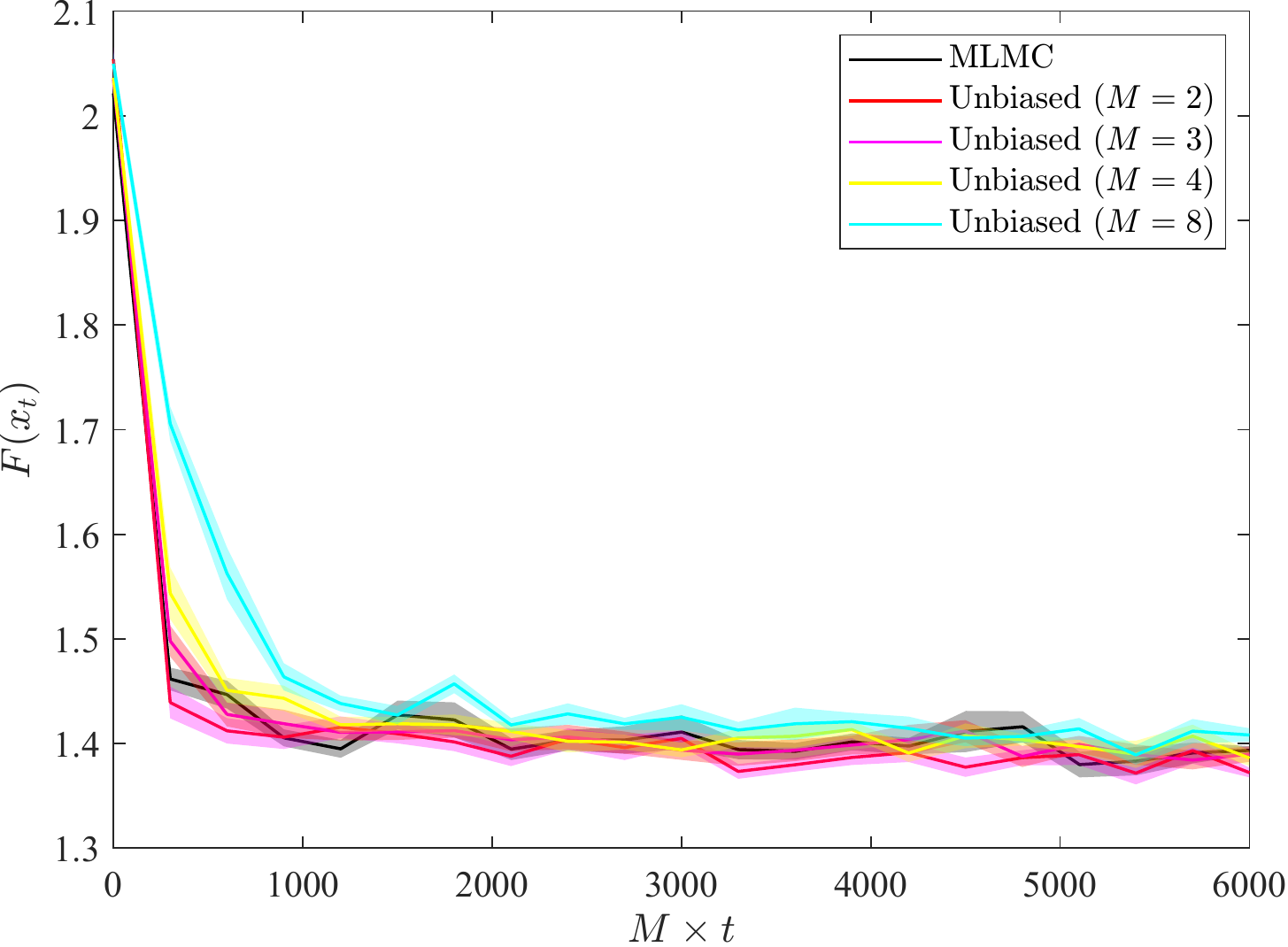}
    \subcaption{$f(x)=|x|$}
    \end{minipage}
    \begin{minipage}{0.45\textwidth}
    \centering
    \includegraphics[width=\textwidth]{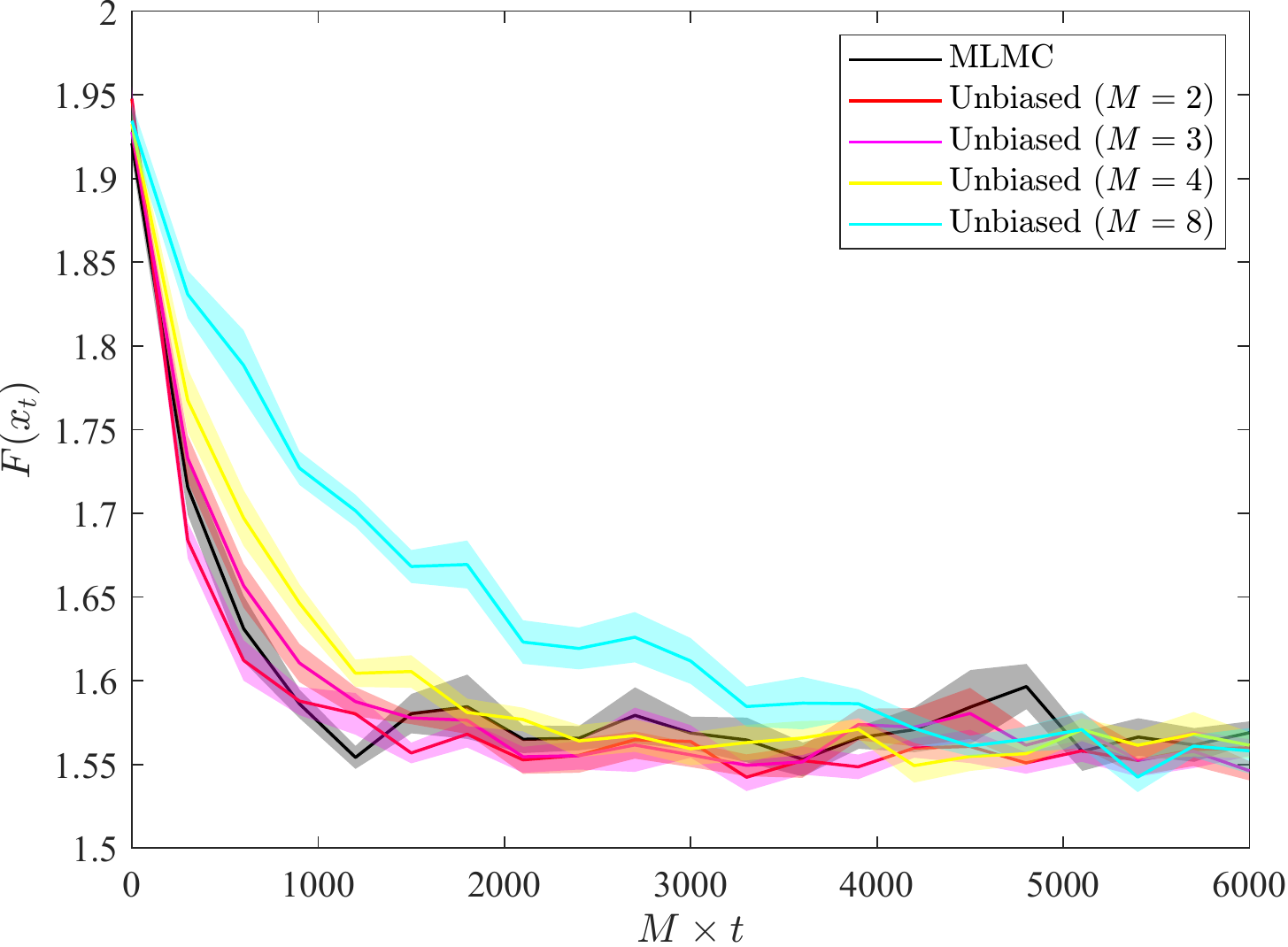}
    \subcaption{$f(x)=\bsone_{x\geq 0}$}
    \end{minipage}
    \caption{\tg{The convergence behavior of the objective function values for the IV regression problem with four different ground truth functions $f$. Several unbiased gradient estimators are used within the Robbins-Monro algorithm. For each gradient estimator, the solid line and the shaded area represent the average and its standard error estimated from 10 independent runs, respectively. Here, Unbiased stands for the third unbiased gradient estimator $\widehat{\nabla F}^{(3)}_{M,N}$.}}
    \label{fig:convergence_iv}
\end{figure}

\tg{Figure~\ref{fig:convergence_iv} shows the behaviors of the objective function value $F(x_t)$. Here the Robbins-Monro iteration with the constant step-sizes equal to $10^{-3}$ terminates after $T$ steps such that the product $MT$ does not exceed $6000$. Since the exact value of $F(x)$ is not available, we estimate it by the unbiased estimator \eqref{eq:special_unbiased} with $N=10^4$ and $M=2$. We draw the solid line to represent the average of 10 estimated values and the shaded area is shown to represent the standard error.}

\tg{Except for $\widehat{\nabla F}^{(3)}_{8,1}$, the difference between the unbiased gradient estimators is not clear for each of the four cases. The convergence of $\widehat{\nabla F}^{(3)}_{8,1}$ is slower, simply because the computational cost per one iteration step is larger than the others. This implies that setting a small value for $M$ is a good choice and saves the overall cost. Although the unbiased MLMC estimator works well at early iteration steps, some minor outliers are detected at later steps. Therefore, using $\widehat{\nabla F}^{(3)}_{M,N}$ should be better for this special case.}

\begin{figure}
    \centering
    \begin{minipage}{0.45\textwidth}
    \centering
    \includegraphics[width=\textwidth]{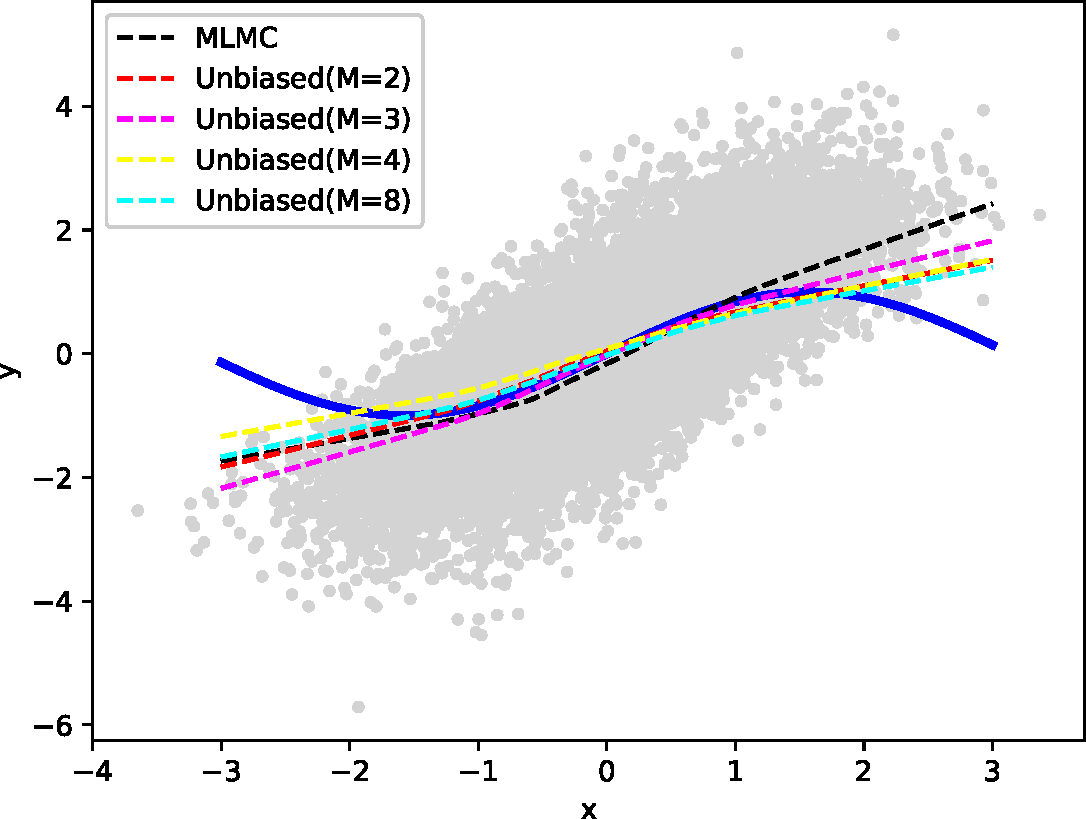}
    \subcaption{$f(x)=\sin(x)$}
    \end{minipage}
    \begin{minipage}{0.45\textwidth}
    \centering
    \includegraphics[width=\textwidth]{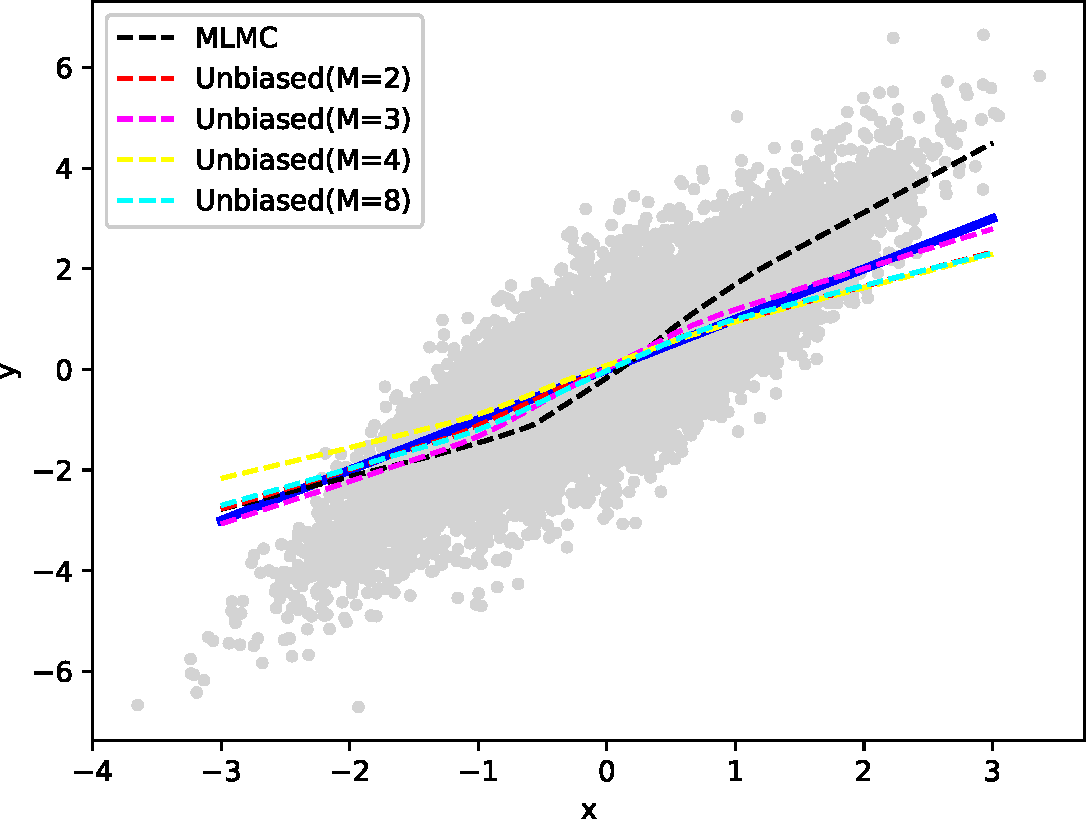}
    \subcaption{$f(x)=x$}
    \end{minipage}\\[5pt]
    \begin{minipage}{0.45\textwidth}
    \centering
    \includegraphics[width=\textwidth]{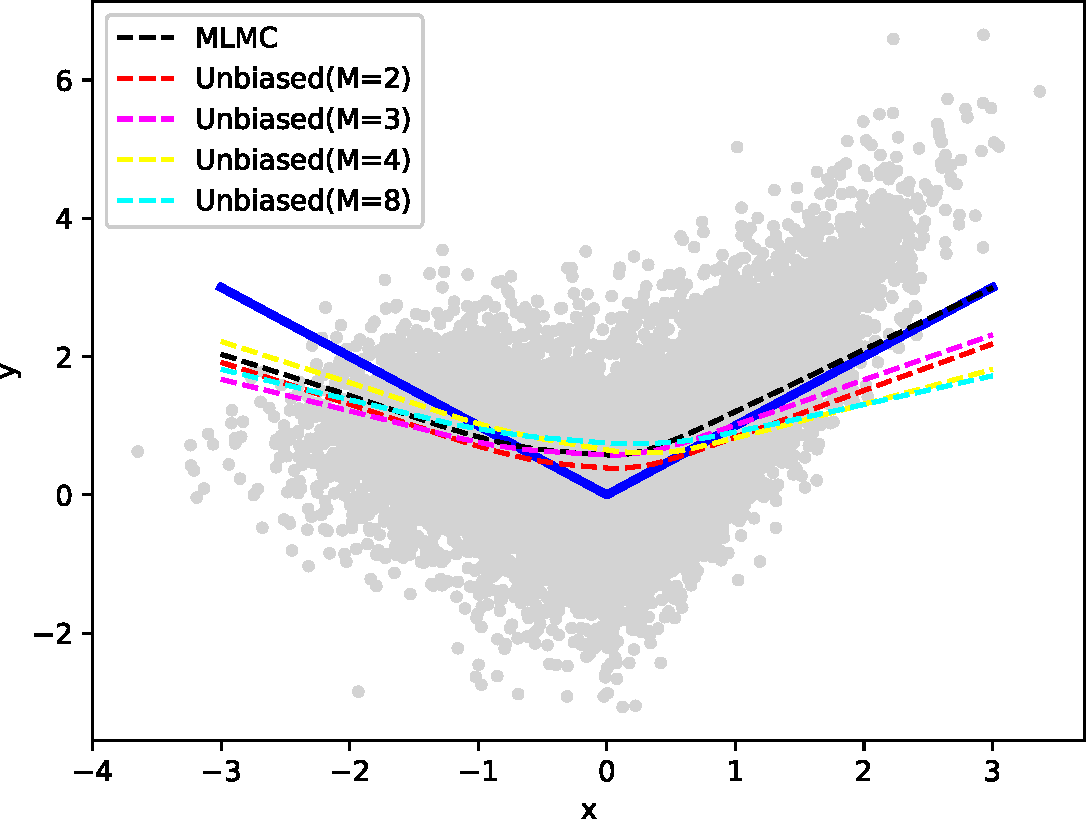}
    \subcaption{$f(x)=|x|$}
    \end{minipage}
    \begin{minipage}{0.45\textwidth}
    \centering
    \includegraphics[width=\textwidth]{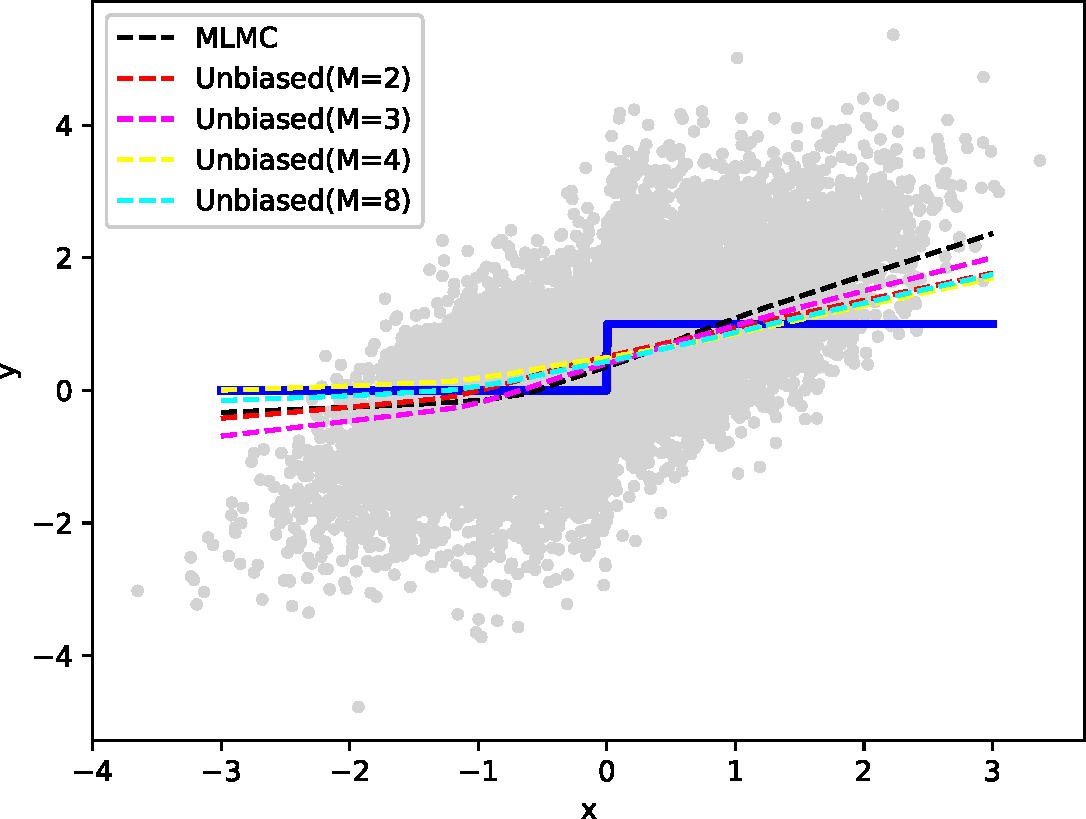}
    \subcaption{$f(x)=\bsone_{x\geq 0}$}
    \end{minipage}
    \caption{\tg{Comparison of the ground truth function $f$ (bold blue) and the neural networks optimized by using several gradient estimators. A set of $10^4$ points for $(X,Y)$ randomly generated by the given process is also plotted. Again, Unbiased stands for the third unbiased gradient estimator $\widehat{\nabla F}^{(3)}_{M,N}$.}}
    \label{fig:signal}
\end{figure}

\tg{Finally we compare the ground truth function $f$ and the optimized neural networks $g_x$ found by the different unbiased gradient estimators in Figure~\ref{fig:signal}. We also plot a set of randomly generated $10^4$ samples for $(X,Y)$ as a reference. Here again, the difference between the gradient estimators is not so clear, except for the case $f(x)=x$ in which the neural network found by the MLMC estimator gives a slightly larger value of $y$ for the region $x>0$. The results for the case $f(x)=\bsone_{x\geq 0}$ are not discontinuous but capture the overall tendency well. Improving the performance for the IV regression itself is an interesting research topic, but of course, is beyond the scope of this paper.}

\section{Concluding remarks}
In this paper, we have studied an unbiased gradient estimation for solving conditional stochastic optimization problems by stochastic gradient descent. In a general problem setting, we have shown the sufficient conditions under which the antithetic property of our unbiased MLMC gradient estimator can be exploited well so that both the variance and the expected computational cost of the estimator can be bounded. Numerical experiments for a simple example confirm a superiority of our MLMC estimator over the standard, nested Monte Carlo estimator. Moreover we have introduced three different unbiased gradient estimators, specific to an important case where the outer function measures the squared loss between an observable variable and a predictor involving an inner conditional expectation. We have shown that the third estimator has minimal variance under a given computational cost.

We leave the following problems open for future research:
\begin{itemize}
    \item Regarding the unbiased MLMC estimator, it might be possible to weaken the sufficient conditions in various ways for more specific problem settings. It is also interesting to investigate what are the necessary conditions for the estimator to have finite variance and expected cost.
    \item We have implicitly assumed that an i.i.d.\ sampling of the inner variable $\eta$ conditioned on the outer variable $\xi$ is possible. In some applications, however, this is not necessarily the case. It is worth studying the use of (unbiased) Markov chain Monte Carlo sampling within our unbiased gradient estimators. A recent work of \cite{WW22} seems quite relevant to this research direction.
\end{itemize}

\section*{Acknowledgements}
The authors would like to thank Yifan Hu (University of Illinois at Urbana-Champaign) and Kei Ishikawa (ETH Z\"{u}rich) for useful discussions and comments. \tg{The authors are also grateful to the anonymous reviewers for their valuable comments} \tgg{and for pointing out an error in the original proof of Theorem~\ref{thm:main1}.}

\appendix
\def\thesection{Appendix~\Alph{section}}

\section{Proof of Theorem~\ref{thm:main1}}\label{app:proof_main1}
\tg{Regarding our matrix norm, it holds for any matrix $A$ and $1\leq p<q$ that
\[ \|A\|_q\leq \|A\|_p\leq N_A^{1/p-1/q}\|A\|_q,\]
where $N_A$ denotes the number of elements of matrix $A$, which will be often used in the subsequent argument.
}

\tg{Moreover, we shall use the following elementary moment inequality, which states that, for a real-valued random variable $X$ with mean 0 and $\EE[|X|^m]<\infty$ for some $m\geq 2$, there exists a constant $C_m>0$ such that
\begin{align}\label{lem:concentration}
    \EE\left[ \left|\frac{1}{N}\sum_{n=1}^{N}X^{(n)}\right|^m\right]\leq C_m\frac{\EE[|X|^m]}{N^{m/2}}
\end{align}
holds for all $N\in \NN$, where $X^{(1)},\ldots,X^{(N)}$ denote the i.i.d.\ copies of $X$, see, e.g., \cite[Chapter~7]{D04}.}

\begin{proof}[Proof of Theorem~\ref{thm:main1}]
For the notational simplicity, we omit the argument $(x,\xi)$ from $g_{\bullet}$ and $\nabla g_{\bullet}$ thoughout this proof. Using the antithetic property \eqref{eq:antithetic2}, we have
\begin{align*}
\Delta \psi_{\ell} & = \left(\overline{\nabla g_{\bullet}}^{\ell}\right)^{\top}\nabla f_{\xi}\left(\overline{g_{\bullet}}^{\ell}\right)-\frac{1}{2}\left( \left(\overline{\nabla g_{\bullet}}^{\ell-1,(a)}\right)^{\top}\nabla f_{\xi}\left(\overline{g_{\bullet}}^{\ell-1,(a)}\right)+\left(\overline{\nabla g_{\bullet}}^{\ell-1,(b)}\right)^{\top}\nabla f_{\xi}\left(\overline{g_{\bullet}}^{\ell-1,(b)}\right)\right)\\
& = \left(\overline{\nabla g_{\bullet}}^{\ell}-\EE_{\eta\mid \xi}\left[ \nabla g_{\eta}\right]\right)^{\top}\left(\nabla f_{\xi}\left(\overline{g_{\bullet}}^{\ell}\right)-\nabla f_{\xi}\left( \EE_{\eta\mid \xi}\left[ g_{\eta}\right]\right)\right) \\
& \quad - \frac{1}{2}\left(\overline{\nabla g_{\bullet}}^{\ell-1,(a)}-\EE_{\eta\mid \xi}\left[ \nabla g_{\eta}\right]\right)^{\top}\left(\nabla f_{\xi}\left(\overline{g_{\bullet}}^{\ell-1,(a)}\right)-\nabla f_{\xi}\left( \EE_{\eta\mid \xi}\left[ g_{\eta}\right]\right)\right) \\
& \quad - \frac{1}{2}\left(\overline{\nabla g_{\bullet}}^{\ell-1,(b)}-\EE_{\eta\mid \xi}\left[ \nabla g_{\eta}\right]\right)^{\top}\left(\nabla f_{\xi}\left(\overline{g_{\bullet}}^{\ell-1,(b)}\right)-\nabla f_{\xi}\left( \EE_{\eta\mid \xi}\left[ g_{\eta}\right]\right)\right) \\
& \quad + \left(\EE_{\eta\mid \xi}\left[ \nabla g_{\eta}\right]\right)^{\top}\left( \nabla f_{\xi}\left(\overline{g_{\bullet}}^{\ell}\right)-\frac{1}{2}\nabla f_{\xi}\left(\overline{g_{\bullet}}^{\ell-1,(a)}\right)-\frac{1}{2}\nabla f_{\xi}\left(\overline{g_{\bullet}}^{\ell-1,(b)}\right)\right).
\end{align*}
Applying Jensen's inequality, the Cauchy–Schwarz inequality, the elementary inequality that $2ab\leq a^2+b^2$ for any $a,b\in \RR$ and the assumption \eqref{assump1} in this order, we get
\begin{align}
\left\|\Delta \psi_{\ell}\right\|_2^2 & \leq 3 \left\|\left(\overline{\nabla g_{\bullet}}^{\ell}-\EE_{\eta\mid \xi}\left[ \nabla g_{\eta}\right]\right)^{\top}\left(\nabla f_{\xi}\left(\overline{g_{\bullet}}^{\ell}\right)-\nabla f_{\xi}\left( \EE_{\eta\mid \xi}\left[ g_{\eta}\right]\right)\right)\right\|_2^2 \notag \\
& \quad + \frac{3}{2}\left\|\left(\overline{\nabla g_{\bullet}}^{\ell-1,(a)}-\EE_{\eta\mid \xi}\left[ \nabla g_{\eta}\right]\right)^{\top}\left(\nabla f_{\xi}\left(\overline{g_{\bullet}}^{\ell-1,(a)}\right)-\nabla f_{\xi}\left( \EE_{\eta\mid \xi}\left[ g_{\eta}\right]\right)\right)\right\|_2^2 \notag \\
& \quad + \frac{3}{2}\left\|\left(\overline{\nabla g_{\bullet}}^{\ell-1,(b)}-\EE_{\eta\mid \xi}\left[ \nabla g_{\eta}\right]\right)^{\top}\left(\nabla f_{\xi}\left(\overline{g_{\bullet}}^{\ell-1,(b)}\right)-\nabla f_{\xi}\left( \EE_{\eta\mid \xi}\left[ g_{\eta}\right]\right)\right)\right\|_2^2 \notag \\
& \quad + 3\left\|\left(\EE_{\eta\mid \xi}\left[ \nabla g_{\eta}\right]\right)^{\top}\left( \nabla f_{\xi}\left(\overline{g_{\bullet}}^{\ell}\right)-\frac{1}{2}\nabla f_{\xi}\left(\overline{g_{\bullet}}^{\ell-1,(a)}\right)-\frac{1}{2}\nabla f_{\xi}\left(\overline{g_{\bullet}}^{\ell-1,(b)}\right)\right)\right\|_2^2 \notag \\
& \leq 3 \left\|\overline{\nabla g_{\bullet}}^{\ell}-\EE_{\eta\mid \xi}\left[ \nabla g_{\eta}\right]\right\|_2^2\left\|\nabla f_{\xi}\left(\overline{g_{\bullet}}^{\ell}\right)-\nabla f_{\xi}\left( \EE_{\eta\mid \xi}\left[ g_{\eta}\right]\right)\right\|_2^2 \notag \\
& \quad + \frac{3}{2}\left\|\overline{\nabla g_{\bullet}}^{\ell-1,(a)}-\EE_{\eta\mid \xi}\left[ \nabla g_{\eta}\right]\right\|_2^2\left\|\nabla f_{\xi}\left(\overline{g_{\bullet}}^{\ell-1,(a)}\right)-\nabla f_{\xi}\left( \EE_{\eta\mid \xi}\left[ g_{\eta}\right]\right)\right\|_2^2 \notag \\
& \quad + \frac{3}{2}\left\|\overline{\nabla g_{\bullet}}^{\ell-1,(b)}-\EE_{\eta\mid \xi}\left[ \nabla g_{\eta}\right]\right\|_2^2\left\|\nabla f_{\xi}\left(\overline{g_{\bullet}}^{\ell-1,(b)}\right)-\nabla f_{\xi}\left( \EE_{\eta\mid \xi}\left[ g_{\eta}\right]\right)\right\|_2^2 \notag \\
& \quad + 3\left\|\EE_{\eta\mid \xi}\left[ \nabla g_{\eta}\right]\right\|_2^2\left\| \nabla f_{\xi}\left(\overline{g_{\bullet}}^{\ell}\right)-\frac{1}{2}\nabla f_{\xi}\left(\overline{g_{\bullet}}^{\ell-1,(a)}\right)-\frac{1}{2}\nabla f_{\xi}\left(\overline{g_{\bullet}}^{\ell-1,(b)}\right)\right\|_2^2 \notag \\
& \leq \frac{3}{2} \left\|\overline{\nabla g_{\bullet}}^{\ell}-\EE_{\eta\mid \xi}\left[ \nabla g_{\eta}\right]\right\|_2^4+\frac{3}{2}\left\|\nabla f_{\xi}\left(\overline{g_{\bullet}}^{\ell}\right)-\nabla f_{\xi}\left( \EE_{\eta\mid \xi}\left[ g_{\eta}\right]\right)\right\|_2^4 \notag \\
& \quad + \frac{3}{4}\left\|\overline{\nabla g_{\bullet}}^{\ell-1,(a)}-\EE_{\eta\mid \xi}\left[ \nabla g_{\eta}\right]\right\|_2^4+\frac{3}{4}\left\|\nabla f_{\xi}\left(\overline{g_{\bullet}}^{\ell-1,(a)}\right)-\nabla f_{\xi}\left( \EE_{\eta\mid \xi}\left[ g_{\eta}\right]\right)\right\|_2^4 \notag \\
& \quad + \frac{3}{4}\left\|\overline{\nabla g_{\bullet}}^{\ell-1,(b)}-\EE_{\eta\mid \xi}\left[ \nabla g_{\eta}\right]\right\|_2^4+\frac{3}{4}\left\|\nabla f_{\xi}\left(\overline{g_{\bullet}}^{\ell-1,(b)}\right)-\nabla f_{\xi}\left( \EE_{\eta\mid \xi}\left[ g_{\eta}\right]\right)\right\|_2^4 \notag \\
& \quad + 3\left\|\EE_{\eta\mid \xi}\left[ \nabla g_{\eta}\right]\right\|_2^2\left\| \nabla f_{\xi}\left(\overline{g_{\bullet}}^{\ell}\right)-\frac{1}{2}\nabla f_{\xi}\left(\overline{g_{\bullet}}^{\ell-1,(a)}\right)-\frac{1}{2}\nabla f_{\xi}\left(\overline{g_{\bullet}}^{\ell-1,(b)}\right)\right\|_2^2 \notag \\
& \leq \frac{3}{2}\left[\left\|\overline{\nabla g_{\bullet}}^{\ell}-\EE_{\eta\mid \xi}\left[ \nabla g_{\eta}\right]\right\|_2^4+\frac{1}{2}\left\|\overline{\nabla g_{\bullet}}^{\ell-1,(a)}-\EE_{\eta\mid \xi}\left[ \nabla g_{\eta}\right]\right\|_2^4+\frac{1}{2}\left\|\overline{\nabla g_{\bullet}}^{\ell-1,(b)}-\EE_{\eta\mid \xi}\left[ \nabla g_{\eta}\right]\right\|_2^4\right] \notag \\
& \quad + \frac{3}{2}\tg{\lambda_1^4} \left[\left\|\overline{g_{\bullet}}^{\ell}-\EE_{\eta\mid \xi}\left[g_{\eta}\right]\right\|_2^4+\frac{1}{2}\left\|\overline{g_{\bullet}}^{\ell-1,(a)}-\EE_{\eta\mid \xi}\left[g_{\eta}\right]\right\|_2^4+\frac{1}{2}\left\|\overline{g_{\bullet}}^{\ell-1,(b)}-\EE_{\eta\mid \xi}\left[g_{\eta}\right]\right\|_2^4\right] \notag \\
& \quad + 3\left\|\EE_{\eta\mid \xi}\left[ \nabla g_{\eta}\right]\right\|_2^2 \left\| \nabla f_{\xi}\left(\overline{g_{\bullet}}^{\ell}\right)-\frac{1}{2}\nabla f_{\xi}\left(\overline{g_{\bullet}}^{\ell-1,(a)}\right)-\frac{1}{2}\nabla f_{\xi}\left(\overline{g_{\bullet}}^{\ell-1,(b)}\right)\right\|_2^2.\label{eq:bound1}
\end{align}

Let us recall that $x\in \Xcal \subseteq \RR^d$ and $g_{\eta}(\cdot,\xi)$ is a function from $\RR^d$ to $\RR^k$.
For the terms in the first and second parentheses of \eqref{eq:bound1}, by using Jensen's inequality, the inequality~\tgg{\eqref{lem:concentration}} with $m=4$ and the tower property of conditional expectation, it holds that
\begin{align*}
    \EE\left[\left\|\overline{\nabla g_{\bullet}}^{\ell}-\EE_{\eta\mid \xi}\left[ \nabla g_{\eta}\right]\right\|_2^4 \right] & = \EE\left[\EE\left[\left\|\overline{\nabla g_{\bullet}}^{\ell}-\EE_{\eta\mid \xi}\left[ \nabla g_{\eta}\right]\right\|_2^4 \mid \xi\right]\right] \\
    & \leq dk\EE\left[\EE\left[\left\|\overline{\nabla g_{\bullet}}^{\ell}-\EE_{\eta\mid \xi}\left[ \nabla g_{\eta}\right]\right\|_4^4 \mid \xi\right]\right] \leq \frac{dkC_4}{2^{2\ell}}\EE_{\xi}\EE_{\eta\mid \xi}\left[\left\|\nabla g_{\eta}-\EE_{\eta\mid \xi}\left[ \nabla g_{\eta}\right]\right\|_4^4 \right],
\end{align*}
and 
\begin{align*}
    \EE\left[\left\|\overline{g_{\bullet}}^{\ell}-\EE_{\eta\mid \xi}\left[g_{\eta}\right]\right\|_2^4 \right] & = \EE\left[\EE\left[\left\|\overline{g_{\bullet}}^{\ell}-\EE_{\eta\mid \xi}\left[g_{\eta}\right]\right\|_2^4 \mid \xi\right]\right] \\
    & \leq k\EE\left[\EE\left[\left\|\overline{g_{\bullet}}^{\ell}-\EE_{\eta\mid \xi}\left[g_{\eta}\right]\right\|_4^4 \mid \xi\right]\right] \leq \frac{kC_4}{2^{2\ell}}\EE_{\xi}\EE_{\eta\mid \xi}\left[\left\|g_{\eta}-\EE_{\eta\mid \xi}\left[g_{\eta}\right]\right\|_4^4 \right].
\end{align*}
Similar bounds also hold for $\overline{\nabla g_{\bullet}}^{\ell-1,(a)}, \overline{\nabla g_{\bullet}}^{\ell-1,(b)}, \overline{g_{\bullet}}^{\ell-1,(a)}$ and $\overline{g_{\bullet}}^{\ell-1,(b)}$. Thus it suffices to show a bound on the expectation of the last term of \eqref{eq:bound1}.

\tgg{For each coordinate index $i$ with unit vector $e_i$, by} applying the mean value theorem for functions with multiple variables and then using the antithetic property \eqref{eq:antithetic1}, \tgg{we see that there exist $t_i^{(a)},t_i^{(b)}\in (0,1)$} such that
\begin{align*}
    & \tgg{e_i}\cdot \left(\nabla f_{\xi}\left(\overline{g_{\bullet}}^{\ell}\right)-\frac{1}{2}\nabla f_{\xi}\left(\overline{g_{\bullet}}^{\ell-1,(a)}\right)-\frac{1}{2}\nabla f_{\xi}\left(\overline{g_{\bullet}}^{\ell-1,(b)}\right)\right) \\
    & = \frac{\tgg{e_i}}{2}\cdot \nabla^{\top}\nabla f_{\xi}\left( \overline{g_{\bullet}}^{\ell-1,(a)}+\tgg{t_i^{(a)}}\left(\overline{g_{\bullet}}^{\ell}-\overline{g_{\bullet}}^{\ell-1,(a)} \right)\right)\cdot\left(\overline{g_{\bullet}}^{\ell}-\overline{g_{\bullet}}^{\ell-1,(a)}\right) \\
    & \quad + \frac{\tgg{e_i}}{2}\nabla^{\top}\nabla f_{\xi}\left( \overline{g_{\bullet}}^{\ell-1,(b)}+\tgg{t_i^{(b)}}\left(\overline{g_{\bullet}}^{\ell}-\overline{g_{\bullet}}^{\ell-1,(b)} \right)\right)\cdot\left(\overline{g_{\bullet}}^{\ell}-\overline{g_{\bullet}}^{\ell-1,(b)}\right) \\
    & = \frac{\tgg{e_i}}{4}\cdot \left(\nabla^{\top}\nabla f_{\xi}\left( \overline{g_{\bullet}}^{\ell-1,(a)}+\tgg{t_i^{(a)}} \left(\overline{g_{\bullet}}^{\ell}-\overline{g_{\bullet}}^{\ell-1,(a)} \right)\right)-\nabla^{\top}\nabla f_{\xi}\left( \overline{g_{\bullet}}^{\ell-1,(b)}+\tgg{t_i^{(b)}} \left(\overline{g_{\bullet}}^{\ell}-\overline{g_{\bullet}}^{\ell-1,(b)} \right)\right) \right)\\
    & \quad \quad \cdot \left(\overline{g_{\bullet}}^{\ell-1,(b)}-\overline{g_{\bullet}}^{\ell-1,(a)} \right).
\end{align*}
The Cauchy–Schwarz inequality and the assumption \eqref{assump2} lead to
\begin{align*}
    & \left|\tgg{e_i}\cdot \left(\nabla f_{\xi}\left(\overline{g_{\bullet}}^{\ell}\right)-\frac{1}{2}\nabla f_{\xi}\left(\overline{g_{\bullet}}^{\ell-1,(a)}\right)-\frac{1}{2}\nabla f_{\xi}\left(\overline{g_{\bullet}}^{\ell-1,(b)}\right)\right)\right|_2^2 \\
    & \leq \frac{1}{16}\left\|\nabla^{\top}\nabla f_{\xi}\left( \overline{g_{\bullet}}^{\ell-1,(a)}+\tgg{t_i^{(a)}} \left(\overline{g_{\bullet}}^{\ell}-\overline{g_{\bullet}}^{\ell-1,(a)} \right)\right)-\nabla^{\top}\nabla f_{\xi}\left( \overline{g_{\bullet}}^{\ell-1,(b)}+\tgg{t_i^{(b)}} \left(\overline{g_{\bullet}}^{\ell}-\overline{g_{\bullet}}^{\ell-1,(b)} \right)\right) \right\|_2^2\\
    & \quad \quad \cdot \left\|\overline{g_{\bullet}}^{\ell-1,(b)}-\overline{g_{\bullet}}^{\ell-1,(a)} \right\|_2^2\\
    & \leq \frac{\lambda_2^2}{16}\left\| \overline{g_{\bullet}}^{\ell-1,(a)}+\tgg{t_i^{(a)}} \left(\overline{g_{\bullet}}^{\ell}-\overline{g_{\bullet}}^{\ell-1,(a)} \right)- \overline{g_{\bullet}}^{\ell-1,(b)}-\tgg{t_i^{(b)}} \left(\overline{g_{\bullet}}^{\ell}-\overline{g_{\bullet}}^{\ell-1,(b)} \right) \right\|_2^{2\rho}\\
    & \quad \quad \cdot \left\|\overline{g_{\bullet}}^{\ell-1,(b)}-\overline{g_{\bullet}}^{\ell-1,(a)} \right\|_2^2 \\
    & \leq \frac{\lambda_2^2}{16} \tgg{\left| \frac{t_i^{(a)}+t_i^{(b)}}{2}-1\right|^{2\rho}}\left\|\overline{g_{\bullet}}^{\ell-1,(b)}-\overline{g_{\bullet}}^{\ell-1,(a)} \right\|_2^{2(1+\rho)} \\
    & \leq \frac{\lambda_2^2}{16}\left\|\overline{g_{\bullet}}^{\ell-1,(b)}-\overline{g_{\bullet}}^{\ell-1,(a)} \right\|_2^{2(1+\rho)}.
\end{align*}
\tgg{Noting that the last bound is independent of the coordinate index $i$, we obtain
\begin{align*}
    & \left\|\nabla f_{\xi}\left(\overline{g_{\bullet}}^{\ell}\right)-\frac{1}{2}\nabla f_{\xi}\left(\overline{g_{\bullet}}^{\ell-1,(a)}\right)-\frac{1}{2}\nabla f_{\xi}\left(\overline{g_{\bullet}}^{\ell-1,(b)}\right)\right\|_2^2 \\
    & \leq \frac{\lambda_2^2 k}{16}\left\|\overline{g_{\bullet}}^{\ell-1,(b)}-\overline{g_{\bullet}}^{\ell-1,(a)} \right\|_2^{2(1+\rho)} \\
    & \leq \frac{\lambda_2^2 k^{1+\rho}}{16}\left\|\overline{g_{\bullet}}^{\ell-1,(b)}-\overline{g_{\bullet}}^{\ell-1,(a)} \right\|_{2(1+\rho)}^{2(1+\rho)}\\
    & \leq \lambda_2^2 k^{1+\rho}2^{2(1+\rho)-5}\left[\left\|\overline{g_{\bullet}}^{\ell-1,(a)}-\EE_{\eta\mid \xi}\left[g_{\eta}\right] \right\|_{2(1+\rho)}^{2(1+\rho)}+\left\|\overline{g_{\bullet}}^{\ell-1,(b)}-\EE_{\eta\mid \xi}\left[g_{\eta}\right] \right\|_{2(1+\rho)}^{2(1+\rho)}\right].
\end{align*}}

Now, for the last term of \eqref{eq:bound1}, we use H\"{o}lder's inequality and the Cauchy–Schwarz inequality, and then apply the inequality~\tgg{\eqref{lem:concentration}} with $m=4(1+\rho)$ to get
\begin{align*}
    & \EE\left[\left\|\EE_{\eta\mid \xi}\left[ \nabla g_{\eta}\right]\right\|_2^2 \left\| \nabla f_{\xi}\left(\overline{g_{\bullet}}^{\ell}\right)-\frac{1}{2}\nabla f_{\xi}\left(\overline{g_{\bullet}}^{\ell-1,(a)}\right)-\frac{1}{2}\nabla f_{\xi}\left(\overline{g_{\bullet}}^{\ell-1,(b)}\right)\right\|_2^2\right] \\
    & \leq \lambda_2^2 \tgg{k^{1+\rho}}2^{2(1+\rho)-5}\EE\left[\left\|\EE_{\eta\mid \xi}\left[ \nabla g_{\eta}\right]\right\|_2^2 \left[\left\|\overline{g_{\bullet}}^{\ell-1,(a)}-\EE_{\eta\mid \xi}\left[g_{\eta}\right] \right\|_{2(1+\rho)}^{2(1+\rho)}+\left\|\overline{g_{\bullet}}^{\ell-1,(b)}-\EE_{\eta\mid \xi}\left[g_{\eta}\right] \right\|_{2(1+\rho)}^{2(1+\rho)}\right]\right] \\
    & \leq \lambda_2^2 \tgg{k^{1+\rho}}2^{2(1+\rho)-5}\left(\EE\left[\left\|\EE_{\eta\mid \xi}\left[ \nabla g_{\eta}\right]\right\|_2^{4}\right]\right)^{1/2} \\
    & \quad \times \left(\left(\EE\left[\left\|\overline{g_{\bullet}}^{\ell-1,(a)}-\EE_{\eta\mid \xi}\left[g_{\eta}\right] \right\|_{2(1+\rho)}^{4(1+\rho)}\right]\right)^{1/2}+\left(\EE\left[\left\|\overline{g_{\bullet}}^{\ell-1,(b)}-\EE_{\eta\mid \xi}\left[g_{\eta}\right] \right\|_{2(1+\rho)}^{4(1+\rho)}\right]\right)^{1/2}\right) \\
    & \leq \lambda_2^2 d^{1/2} \tgg{k^{2+\rho}}2^{2(1+\rho)-5}\left(\EE\left[\left\|\EE_{\eta\mid \xi}\left[ \nabla g_{\eta}\right]\right\|_{4}^{4}\right]\right)^{1/2} \\
    & \quad \times \left(\left(\EE\left[\left\|\overline{g_{\bullet}}^{\ell-1,(a)}-\EE_{\eta\mid \xi}\left[g_{\eta}\right] \right\|_{4(1+\rho)}^{4(1+\rho)}\right]\right)^{1/2}+\left(\EE\left[\left\|\overline{g_{\bullet}}^{\ell-1,(b)}-\EE_{\eta\mid \xi}\left[g_{\eta}\right] \right\|_{4(1+\rho)}^{4(1+\rho)}\right]\right)^{1/2}\right) \\
    & \leq \frac{\lambda_2^2 d^{1/2} \tgg{k^{2+\rho}}2^{2(1+\rho)-4}C_{4(1+\rho)}^{1/2}}{2^{(1+\rho)(\ell-1)}}\left(\EE_{\xi}\EE_{\eta\mid \xi}\left[\left\|\nabla g_{\eta}\right\|_{4}^{4}\right]\right)^{1/2}\left(\EE_{\xi}\EE_{\eta\mid \xi}\left[\left\|g_{\eta}-\EE_{\eta\mid \xi}\left[ g_{\eta}\right]\right\|_{4(1+\rho)}^{4(1+\rho)} \right]\right)^{1/2}.
\end{align*}
The last upper bound is finite as we assume $p\geq 4$ and $q\geq 4(1+\rho)$. \tg{Here we note that, for the second inequality, we have used the simplest version of H\"{o}lder's inequality, i.e., H\"{o}lder's inequality with H\"{o}lder conjugates $p=q=2$. One can replace it with a more general version of H\"{o}lder's inequality and still obtain a similar upper bound which tightens the restriction on $p$ and loosens the restriction on $q$, or vice versa.}

Altogether we obtain a bound on the expected squared 2-norm of $\Delta \psi_{\ell}(x)$ as
\begin{align*}
    & \EE\left[\left\|\Delta \psi_{\ell}(x)\right\|_2^2\right] \\
    & \leq \frac{3}{2}dkC_4\EE_{\xi}\EE_{\eta\mid \xi}\left[\left\|\nabla g_{\eta}-\EE_{\eta\mid \xi}\left[ \nabla g_{\eta}\right]\right\|_4^4 \right]\left[\frac{1}{2^{2\ell}}+\frac{1}{2^{2(\ell-1)+1}}+\frac{1}{2^{2(\ell-1)+1}}\right] \notag \\
    & \quad + \frac{3}{2}\tg{\lambda_1^4}kC_4\EE_{\xi}\EE_{\eta\mid \xi}\left[\left\|g_{\eta}-\EE_{\eta\mid \xi}\left[g_{\eta}\right]\right\|_4^4 \right] \left[\frac{1}{2^{2\ell}}+\frac{1}{2^{2(\ell-1)+1}}+\frac{1}{2^{2(\ell-1)+1}}\right] \notag \\
    & \quad + 3\frac{\lambda_2^2 d^{1/2} \tgg{k^{2+\rho}}2^{2(1+\rho)-4}C_{4(1+\rho)}^{1/2}}{2^{(1+\rho)(\ell-1)}}\left(\EE_{\xi}\EE_{\eta\mid \xi}\left[\left\|\nabla g_{\eta}\right\|_{4}^{4}\right]\right)^{1/2}\left(\EE_{\xi}\EE_{\eta\mid \xi}\left[\left\|g_{\eta}-\EE_{\eta\mid \xi}\left[ g_{\eta}\right]\right\|_{4(1+\rho)}^{4(1+\rho)} \right]\right)^{1/2}.
\end{align*}
This implies the existence of a constant $c_2>0$, independent of $\ell$, such that $\EE\left[\left\|\Delta \psi_{\ell}(x)\right\|_2^2\right]\leq c_2 2^{-(1+\rho)\ell}$ for all $\ell\in \NN$, which completes the proof.
\end{proof}

\section{Proof of Theorem~\ref{thm:main2}}\label{app:proof_main2}
Since it is obvious from the independence between $\xi^{(1)},\ldots,\xi^{(N)}$ that
\[ \EE\left[ \left\|\widehat{\nabla F}^{(i)}_{M,N}(x)-\nabla F(x)\right\|_2^2\right]=\frac{1}{N}\EE\left[ \left\|\widehat{\nabla F}^{(i)}_{M,1}(x)-\nabla F(x)\right\|_2^2\right],\]
for all $i=1,2,3$, it suffices to prove the theorem only for the case $N=1$. For the notational simplicity, we shall drop the superscript from $\xi^{(1)}$ and write $\eta^{(m)}$ and $\eta'^{(m)}$ instead of $\eta^{(m,1)}$ and $\eta'^{(m,1)}$, respectively. Moreover, for given $M$, we write
\begin{align*}
    & \overline{g}_{\bullet}^{(a)}=\frac{1}{M}\sum_{m=1}^{M}g_{\eta^{(m)}}(x, \xi), \quad \overline{g}_{\bullet}^{(b)}=\frac{1}{M}\sum_{m=1}^{M}g_{\eta'^{(m)}}(x, \xi),\\
    & \overline{\nabla g}_{\bullet}^{(a)}=\frac{1}{M}\sum_{m=1}^{M}\nabla g_{\eta^{(m)}}(x, \xi), \quad \overline{\nabla g}_{\bullet}^{(b)}=\frac{1}{M}\sum_{m=1}^{M}\nabla g_{\eta'^{(m)}}(x, \xi), \quad \text{and}\\
    & \overline{g\nabla g}_{\bullet}^{(a)}=\frac{1}{M}\sum_{m=1}^{M}g_{\eta^{(m)}}(x, \xi)\nabla g_{\eta^{(m)}}(x, \xi),
\end{align*} 
which simplifies the forms of the estimators into
\begin{align*}
    & \widehat{\nabla F}^{(1)}_{M,1}(x) = -2 \left(u(\xi)-\overline{g}_{\bullet}^{(a)}\right)\overline{\nabla g}_{\bullet}^{(b)},\quad \widehat{\nabla F}^{(2)}_{M,1}(x) = - \left(u(\xi)-\overline{g}_{\bullet}^{(a)}\right)\overline{\nabla g}_{\bullet}^{(b)}- \left(u(\xi)-\overline{g}_{\bullet}^{(b)}\right)\overline{\nabla g}_{\bullet}^{(a)},\quad \text{and}\\
    & \widehat{\nabla F}^{(3)}_{M,1}(x) = -2 \left(u(\xi)\overline{\nabla g}_{\bullet}^{(a)}-\frac{M}{M-1}\overline{g}_{\bullet}^{(a)}\overline{\nabla g}_{\bullet}^{(a)}+\frac{1}{M-1}\overline{g\nabla g}_{\bullet}^{(a)}\right),
\end{align*}
respectively. 

The following lemma is crucial in the proof of the theorem. Since it can be proven by elementary calculations, we omit the proof.

\begin{lemma}\label{lem:many_equalities}
With the notation above, we have
\begin{align*}
    & \EE\left[\overline{g}_{\bullet}^{(a)}\mid \xi\right] = \EE\left[\overline{g}_{\bullet}^{(b)}\mid \xi\right] = \EE_{\eta\mid \xi}\left[ g_{\eta}\right], \\
    & \EE\left[\left(\overline{g}_{\bullet}^{(a)}\right)^2\mid \xi\right] = \EE\left[\left(\overline{g}_{\bullet}^{(b)}\right)^2\mid \xi\right] = \left(\EE_{\eta\mid \xi}\left[ g_{\eta}\right]\right)^2+\frac{\EE_{\eta\mid \xi}\left[(g_{\eta}-\EE_{\eta\mid \xi}\left[ g_{\eta}\right])^2\right]}{M}, \\
    & \EE\left[\overline{\nabla g}_{\bullet}^{(a)}\mid \xi\right] = \EE\left[\overline{\nabla  g}_{\bullet}^{(b)}\mid \xi\right] = \EE_{\eta\mid \xi}\left[\nabla  g_{\eta}\right], \\
    & \EE\left[\left\|\overline{\nabla g}_{\bullet}^{(a)}\right\|_2^2\mid \xi\right] = \EE\left[\left\|\overline{\nabla g}_{\bullet}^{(b)}\right\|_2^2\mid \xi\right] = \left\|\EE_{\eta\mid \xi}\left[ \nabla g_{\eta}\right]\right\|_2^2+ \frac{\EE_{\eta\mid \xi}\left[\|\nabla g_{\eta}-\EE_{\eta\mid \xi}\left[ \nabla g_{\eta}\right]\|_2^2\right]}{M}, \\
    & \EE\left[\overline{g}_{\bullet}^{(a)}\overline{\nabla g}_{\bullet}^{(a)}\mid \xi\right] = \EE\left[\overline{g}_{\bullet}^{(b)}\overline{\nabla g}_{\bullet}^{(b)}\mid \xi\right] = \frac{M-1}{M}\EE_{\eta\mid \xi}\left[g_{\eta}\right]\EE_{\eta\mid \xi}\left[\nabla g_{\eta}\right]+\frac{1}{M}\EE_{\eta\mid \xi}\left[ g_{\eta}\nabla g_{\eta}\right],\\
    & \EE\left[\overline{g}_{\bullet}^{(a)}\left\|\overline{\nabla g}_{\bullet}^{(a)}\right\|_2^2\mid \xi\right] = \frac{\EE_{\eta\mid \xi}\left[ g_{\eta}\left\|\nabla g_{\eta}\right\|_2^2\right]}{M^2}+\frac{M-1}{M^2}\EE_{\eta\mid \xi}\left[ g_{\eta}\right]\EE_{\eta\mid \xi}\left[\left\|\nabla g_{\eta}\right\|_2^2\right]\\
    & \quad \quad \quad \quad \quad \quad \quad \quad \quad \quad +\frac{2(M-1)}{M^2}\EE_{\eta\mid \xi}\left[g_{\eta}\nabla g_{\eta}\right]\cdot \EE_{\eta\mid \xi}\left[\nabla g_{\eta}\right]\\
    & \quad \quad \quad \quad \quad \quad \quad \quad \quad \quad + \frac{(M-1)(M-2)}{M^2}\EE_{\eta\mid \xi}\left[g_{\eta}\right]\left\|\EE_{\eta\mid \xi}\left[\nabla g_{\eta}(X)\right]\right\|_2^2, \\
    & \EE\left[\left(\overline{g}_{\bullet}^{(a)}\right)^2\left\|\overline{\nabla g}_{\bullet}^{(a)}\right\|_2^2\mid \xi\right] = \frac{\EE_{\eta\mid \xi}\left[\left\|g_{\eta}\nabla g_{\eta}\right\|_2^2\right]}{M^3} + \frac{2(M-1)}{M^3}\EE_{\eta\mid \xi}\left[g_{\eta}\right]\EE_{\eta\mid \xi}\left[g_{\theta}\left\|\nabla g_{\eta}\right\|_2^2\right]\\
    & \quad \quad \quad \quad \quad \quad \quad \quad \quad \quad \quad \quad + \frac{2(M-1)}{M^3}\EE_{\eta\mid \xi}\left[\left(g_{\eta}\right)^2\nabla g_{\eta}\right]\cdot \EE_{\eta\mid \xi}\left[\nabla g_{\eta}\right]\\
    & \quad \quad \quad \quad \quad \quad \quad \quad \quad \quad \quad \quad + \frac{2(M-1)}{M^3}\left\|\EE_{\eta\mid \xi}\left[g_{\theta}\nabla g_{\eta}\right]\right\|_2^2\\
    & \quad \quad \quad \quad \quad \quad \quad \quad \quad \quad \quad \quad + \frac{M-1}{M^3}\EE_{\eta\mid \xi}\left[(g_{\eta})^2\right]\EE_{\eta\mid \xi}\left[\left\|\nabla g_{\eta}\right\|_2^2\right]\\
    & \quad \quad \quad \quad \quad \quad \quad \quad \quad \quad \quad \quad + \frac{(M-1)(M-2)}{M^3}\EE_{\eta\mid \xi}\left[(g_{\eta})^2\right]\left\|\EE_{\eta\mid \xi}\left[\nabla g_{\eta}\right]\right\|_2^2\\
    & \quad \quad \quad \quad \quad \quad \quad \quad \quad \quad \quad \quad + \frac{(M-1)(M-2)}{M^3}\left(\EE_{\eta\mid \xi}\left[g_{\eta}\right]\right)^2\EE_{\eta\mid \xi}\left[\left\|\nabla g_{\eta}\right\|_2^2\right]\\
    & \quad \quad \quad \quad \quad \quad \quad \quad \quad \quad \quad \quad + \frac{4(M-1)(M-2)}{M^3}\EE_{\eta\mid \xi}\left[g_{\eta}\right]\EE_{\eta\mid \xi}\left[\nabla g_{\eta}\right]\cdot \EE_{\eta\mid \xi}\left[g_{\eta}\nabla g_{\eta}\right]\\
    & \quad \quad \quad \quad \quad \quad \quad \quad \quad \quad \quad \quad + \frac{(M-1)(M-2)(M-3)}{M^3}\left(\EE_{\eta\mid \xi}\left[g_{\eta}\right]\right)^2\left\|\EE_{\eta\mid \xi}\left[\nabla g_{\eta}\right]\right\|_2^2, \\
    & \EE\left[\overline{\nabla g}_{\bullet}^{(a)}\cdot \overline{g\nabla g}_{\bullet}^{(a)} \mid \xi\right] = \frac{\EE_{\eta\mid \xi}\left[g_{\eta}\left\|\nabla g_{\eta}\right\|_2^2\right]}{M}+\frac{M-1}{M}\EE_{\eta\mid \xi}\left[\nabla g_{\eta}\right]\cdot \EE_{\eta\mid \xi}\left[g_{\eta}\nabla g_{\eta}\right], \\
    & \EE\left[\overline{g}_{\bullet}^{(a)}\overline{\nabla g}_{\bullet}^{(a)}\cdot \overline{g\nabla g}_{\bullet}^{(a)} \mid \xi\right] = \frac{\EE_{\eta\mid \xi}\left[\left\|g_{\eta}\nabla g_{\eta}\right\|_2^2\right]}{M^2}+\frac{M-1}{M^2}\left\|\EE_{\eta\mid \xi}\left[g_{\eta}\nabla g_{\eta}\right]\right\|_2^2 \\
    & \quad \quad \quad \quad \quad \quad \quad \quad \quad \quad \quad \quad + \frac{M-1}{M^2}\EE_{\eta\mid \xi}\left[\left(g_{\eta}\right)^2\nabla g_{\eta}\right]\cdot \EE_{\eta\mid \xi}\left[\nabla g_{\eta}\right] \\
    & \quad \quad \quad \quad \quad \quad \quad \quad \quad \quad \quad \quad + \frac{M-1}{M^2}\EE_{\eta\mid \xi}\left[g_{\eta}\right]\EE_{\eta\mid \xi}\left[g_{\eta}\left\|\nabla g_{\eta}\right\|_2^2\right]\\
    & \quad \quad \quad \quad \quad \quad \quad \quad \quad \quad \quad \quad + \frac{(M-1)(M-2)}{M^2}\EE_{\eta\mid \xi}\left[g_{\eta}\right]\EE_{\eta\mid \xi}\left[\nabla g_{\eta}\right]\cdot \EE_{\eta\mid \xi}\left[g_{\eta}\nabla g_{\eta}\right],\\
    & \EE\left[\left\|\overline{g\nabla g}_{\bullet}^{(a)}\right\|_2^2\mid \xi\right] = \frac{\EE_{\eta\mid \xi}\left[\left\|g_{\eta}\nabla g_{\eta}\right\|_2^2\right]}{M}+\frac{M-1}{M}\left\|\EE_{\eta\mid \xi}\left[g_{\eta}\nabla g_{\eta}\right]\right\|_2^2.
\end{align*}
\end{lemma}

\begin{proof}[Proof of Theorem~\ref{thm:main2}]
Since the considered gradient estimators are all unbiased, it suffices to prove
\begin{align}\label{eq:statement_to_be_proven}
    \EE\left[ \left\|\widehat{\nabla F}^{(1)}_{M,1}(x)\right\|_2^2\right] \geq \EE\left[ \left\|\widehat{\nabla F}^{(2)}_{M,1}(x)\right\|_2^2\right]\geq \EE\left[ \left\|\widehat{\nabla F}^{(3)}_{2M,1}(x)\right\|_2^2\right].
\end{align}
Let us consider the first estimator $\widehat{\nabla F}^{(1)}_{M,1}(x)$. Because of the conditional independence between $\overline{g}_{\bullet}^{(a)}$ and $\overline{\nabla g}_{\bullet}^{(b)}$, using Lemma~\ref{lem:many_equalities}, we have
\begin{align*}
    \EE\left[ \left\|\widehat{\nabla F}^{(1)}_{M,1}(x)\right\|_2^2\right] 
    & = 4\EE_{\xi}\left[\EE\left[ \left(u(\xi)-\overline{g}_{\bullet}^{(a)}\right)^2\mid \xi\right]\EE\left[\left\|\overline{\nabla g}_{\bullet}^{(b)}\right\|_2^2 \mid \xi \right]\right] \\
    & = 4\EE_{\xi}\left[\left(\left(u(\xi)- \EE_{\eta\mid \xi}\left[ g_{\eta}\right]\right)^2+\frac{\EE_{\eta\mid \xi}\left[(g_{\eta}-\EE_{\eta\mid \xi}\left[ g_{\eta}\right])^2\right]}{M}\right)\right. \\
    & \quad \quad \quad \quad \times \left. \left( \left\|\EE_{\eta\mid \xi}\left[ \nabla g_{\eta}\right]\right\|_2^2+ \frac{\EE_{\eta\mid \xi}\left[\|\nabla g_{\eta}-\EE_{\eta\mid \xi}\left[ \nabla g_{\eta}\right]\|_2^2\right]}{M}\right) \right].
\end{align*}
Similarly, for the second estimator $\widehat{\nabla F}^{(2)}_{M,1}(x)$, we have
\begin{align*}
    & \EE\left[ \left\|\widehat{\nabla F}^{(2)}_{M,1}(x)\right\|_2^2\right] \\
    & = \EE_{\xi}\left[\EE\left[ \left(u(\xi)-\overline{g}_{\bullet}^{(a)}\right)^2\mid \xi\right]\EE\left[\left\|\overline{\nabla g}_{\bullet}^{(b)}\right\|_2^2 \mid \xi \right]\right]+\EE_{\xi}\left[\EE\left[ \left(u(\xi)-\overline{g}_{\bullet}^{(b)}\right)^2\mid \xi\right]\EE\left[\left\|\overline{\nabla g}_{\bullet}^{(a)}\right\|_2^2 \mid \xi \right]\right] \\
    & \quad + 2\EE_{\xi}\left[\EE\left[ \left(u(\xi)-\overline{g}_{\bullet}^{(a)}\right)\overline{\nabla g}_{\bullet}^{(a)}\mid \xi\right]\cdot \EE\left[\left(u(\xi)-\overline{g}_{\bullet}^{(b)}\right)\overline{\nabla g}_{\bullet}^{(b)}\mid \xi\right]\right] \\
    & = 2\EE_{\xi}\left[\left(\left(u(\xi)- \EE_{\eta\mid \xi}\left[ g_{\eta}\right]\right)^2+\frac{\EE_{\eta\mid \xi}\left[(g_{\eta}-\EE_{\eta\mid \xi}\left[ g_{\eta}\right])^2\right]}{M}\right)\right. \\
    & \quad \quad \quad \quad \times \left. \left( \left\|\EE_{\eta\mid \xi}\left[ \nabla g_{\eta}\right]\right\|_2^2+ \frac{\EE_{\eta\mid \xi}\left[\|\nabla g_{\eta}-\EE_{\eta\mid \xi}\left[ \nabla g_{\eta}\right]\|_2^2\right]}{M}\right) \right. \\
    & \quad \quad \quad \quad + \left. \left\|\left(u(\xi)-\EE_{\eta\mid \xi}\left[g_{\eta}\right]\right) \EE_{\eta\mid \xi}\left[ \nabla g_{\eta}\right]+\frac{1}{M}\left(\EE_{\eta\mid \xi}\left[g_{\eta}\right]\EE_{\eta\mid \xi}\left[\nabla g_{\eta}\right]-\EE_{\eta\mid \xi}\left[ g_{\eta}\nabla g_{\eta}\right]\right)\right\|_2^2\right].
\end{align*}

Now let us compare $\widehat{\nabla F}^{(1)}_{M,1}(x)$ and $\widehat{\nabla F}^{(2)}_{M,1}(x)$. Using the above results and noting that
\begin{align*}
    \left\|\EE_{\eta\mid \xi}\left[ g_{\eta}\nabla g_{\eta}\right]-\EE_{\eta\mid \xi}\left[g_{\eta}\right]\EE_{\eta\mid \xi}\left[\nabla g_{\eta}\right]\right\|_2 & = \left\|\EE_{\eta\mid \xi}\left[ \left(g_{\eta}-\EE_{\eta\mid \xi}\left[g_{\eta}\right]\right)\left(\nabla g_{\eta}-\EE_{\eta\mid \xi}\left[\nabla g_{\eta}\right]\right)\right]\right\|_2\\
    & \leq \left(\EE_{\eta\mid \xi}\left[ \left(g_{\eta}-\EE_{\eta\mid \xi}\left[g_{\eta}\right]\right)^2\right]\EE_{\eta\mid \xi}\left[\left\|\nabla g_{\eta}-\EE_{\eta\mid \xi}\left[\nabla g_{\eta}\right]\right\|_2^2\right]\right)^{1/2},
\end{align*}
we can see that
\begin{align*}
    & \frac{1}{2}\EE\left[ \left\|\widehat{\nabla F}^{(1)}_{M,1}(x)\right\|_2^2-\left\|\widehat{\nabla F}^{(2)}_{M,1}(x)\right\|_2^2\mid \xi\right] \\
    & = \left(\left(u(\xi)- \EE_{\eta\mid \xi}\left[ g_{\eta}\right]\right)^2+\frac{\EE_{\eta\mid \xi}\left[(g_{\eta}-\EE_{\eta\mid \xi}\left[ g_{\eta}\right])^2\right]}{M}\right)\left( \left\|\EE_{\eta\mid \xi}\left[ \nabla g_{\eta}\right]\right\|_2^2+ \frac{\EE_{\eta\mid \xi}\left[\|\nabla g_{\eta}-\EE_{\eta\mid \xi}\left[ \nabla g_{\eta}\right]\|_2^2\right]}{M}\right)\\
    & \quad -\left\|\left(u(\xi)-\EE_{\eta\mid \xi}\left[g_{\eta}\right]\right) \EE_{\eta\mid \xi}\left[ \nabla g_{\eta}\right]+\frac{1}{M}\left(\EE_{\eta\mid \xi}\left[g_{\eta}\right]\EE_{\eta\mid \xi}\left[\nabla g_{\eta}\right]-\EE_{\eta\mid \xi}\left[ g_{\eta}\nabla g_{\eta}\right]\right)\right\|_2^2 \\
    & \geq \frac{1}{M}\left( \left|u(\xi)- \EE_{\eta\mid \xi}\left[ g_{\eta}\right]\right| \left(\EE_{\eta\mid \xi}\left[\left\|\nabla g_{\eta}-\EE_{\eta\mid \xi}\left[\nabla g_{\eta}\right]\right\|_2^2\right]\right)^{1/2}-\left(\EE_{\eta\mid \xi}\left[ \left(g_{\eta}-\EE_{\eta\mid \xi}\left[g_{\eta}\right]\right)^2\right]\right)^{1/2}\left\|\EE_{\eta\mid \xi}\left[ \nabla g_{\eta}\right]\right\|_2\right)^2\\
    & \geq 0.
\end{align*}
Since this inequality holds pointwise for $\xi$, we complete the proof of the first inequality of \eqref{eq:statement_to_be_proven}.

Let us move on to the third estimator $\widehat{\nabla F}^{(3)}_{M,1}(x)$. Using Lemma~\ref{lem:many_equalities}, for $M\geq 2$, we have
\begin{align*}
    \EE\left[ \left\|\widehat{\nabla F}^{(3)}_{M,1}(x)\right\|_2^2\right] & = 4\EE_{\xi}\left[(u(\xi))^2\EE\left[\left\|\overline{\nabla g}_{\bullet}^{(a)}\right\|_2^2\mid \xi\right] +\frac{M^2}{(M-1)^2}\EE\left[\left(\overline{g}_{\bullet}^{(a)}\right)^2\left\|\overline{\nabla g}_{\bullet}^{(a)}\right\|_2^2\mid \xi\right] \right. \\
    & \quad \quad \quad \quad \left. +\frac{1}{(M-1)^2}\EE\left[\left\|\overline{g\nabla g}_{\bullet}^{(a)}\right\|_2^2\mid \xi\right] - \frac{2M}{M-1}u(\xi)\EE\left[\overline{g}_{\bullet}^{(a)}\left\|\overline{\nabla g}_{\bullet}^{(a)}\right\|_2^2\mid \xi\right] \right.\\
    & \quad \quad \quad \quad \left.  +\frac{2}{M-1}u(\xi)\EE\left[\overline{\nabla g}_{\bullet}^{(a)}\cdot \overline{g\nabla g}_{\bullet}^{(a)} \mid \xi\right] -\frac{2M}{(M-1)^2}\EE\left[\overline{g}_{\bullet}^{(a)}\overline{\nabla g}_{\bullet}^{(a)}\cdot \overline{g\nabla g}_{\bullet}^{(a)} \mid \xi\right]\right]\\
    & = 4\EE_{\xi}\left[\left(\left(u(\xi)- \EE_{\eta\mid \xi}\left[ g_{\eta}\right]\right)^2+\frac{\EE_{\eta\mid \xi}\left[(g_{\eta}-\EE_{\eta\mid \xi}\left[ g_{\eta}\right])^2\right]}{M}\right)\right. \\
    & \quad \quad \quad \quad \times \left. \left( \left\|\EE_{\eta\mid \xi}\left[ \nabla g_{\eta}\right]\right\|_2^2+ \frac{\EE_{\eta\mid \xi}\left[\|\nabla g_{\eta}-\EE_{\eta\mid \xi}\left[ \nabla g_{\eta}\right]\|_2^2\right]}{M}\right) \right].\\
    & \quad \quad \quad \quad \left. -\frac{2}{M}\left(u(\xi)-\EE_{\eta\mid \xi}\left[ g_{\eta}\right]\right) \left(  \EE_{\eta\mid \xi}\left[\nabla  g_{\eta}\right]\cdot \EE_{\eta\mid \xi}\left[ g_{\eta}\nabla g_{\eta}\right]-\EE_{\eta\mid \xi}\left[ g_{\eta}\right] \left\|\EE_{\eta\mid \xi}\left[ \nabla g_{\eta}\right]\right\|_2^2\right)\right.\\
    & \quad \quad \quad \quad \left. +\frac{1}{M(M-1)} \left\| \EE_{\eta\mid \xi}\left[ g_{\eta}\right] \EE_{\eta\mid \xi}\left[ \nabla g_{\eta}\right]- \EE_{\eta\mid \xi}\left[ g_{\eta}\nabla g_{\eta}\right]\right\|_2^2\right. \\
    & \quad \quad \quad \quad \left. +\frac{1}{M^2(M-1)} \EE_{\eta\mid \xi}\left[(g_{\eta}-\EE_{\eta\mid \xi}\left[ g_{\eta}\right])^2\right]\EE_{\eta\mid \xi}\left[\|\nabla g_{\eta}-\EE_{\eta\mid \xi}\left[ \nabla g_{\eta}\right]\|_2^2\right]\right].
\end{align*}
For given $\xi$ and any $M\geq 1$, we can see that
\begin{align*}
    & \frac{1}{2}\EE\left[ \left\|\widehat{\nabla F}^{(2)}_{M,1}(x)\right\|_2^2-\left\|\widehat{\nabla F}^{(3)}_{2M,1}(x)\right\|_2^2\mid \xi\right] \\
    & = \frac{M-1}{M^2(2M-1)}\left[\EE_{\eta\mid \xi}\left[(g_{\eta}-\EE_{\eta\mid \xi}\left[ g_{\eta}\right])^2\right]\EE_{\eta\mid \xi}\left[\|\nabla g_{\eta}-\EE_{\eta\mid \xi}\left[ \nabla g_{\eta}\right]\|_2^2\right]+ \left\| \EE_{\eta\mid \xi}\left[ g_{\eta}\right] \EE_{\eta\mid \xi}\left[ \nabla g_{\eta}\right]- \EE_{\eta\mid \xi}\left[ g_{\eta}\nabla g_{\eta}\right]\right\|_2^2\right]\\
    & \geq 0.
\end{align*}
Taking the expectation of both sides with respect to $\xi$ directly completes the proof of the second inequality of \eqref{eq:statement_to_be_proven}. Thus we are done.
\end{proof}


\bibliographystyle{plain}
\bibliography{ref.bib}

\end{document}